\documentstyle[12pt]{article}
\newcommand{\rbox}[1]{{\rm #1}}        
\newcommand{\bbox}[1]{{\bf #1}}
\newcommand{\bsym}[1]{{\bf #1}}   
\newcommand{\rmt}{\rm}
\newcommand{\itt}{\it}
\newcommand{\bft}{\bf}
\newcommand{\relf}[1]{{\sf #1}}
\def\restriction{\mathbin{\hspace{0.1ex}|\hspace{0.1ex}}}
\def\subsetneqq{\mathbin{\hspace{-0.8mm}
\raisebox{-0.8ex}{
$\msur
\stackrel
{\protect\textstyle\subset}
{\scriptstyle\not=}
\msur$
}\hspace{-0.8mm} }}
\newcommand{\bbb}{\hspace{0pt}}
\newcommand{\pri}{{{\bbb\rbox{pr}\bbb}_1\hspace{1.5pt}}}
\newcommand{\prii}{{{\bbb\rbox{pr}\bbb}_2\hspace{1.5pt}}}

\newtheorem{theorem}{Theorem}
\newtheorem{assertion}[theorem]{Assertion}
\newtheorem{corollary}[theorem]{Corollary}
\newtheorem{definition}[theorem]{Definition}
\newtheorem{lemma}[theorem]{Lemma}
\newtheorem{proposition}[theorem]{Proposition}
\newtheorem{remark}[theorem]{Remark}
\newcommand{\proof}{\noi{\bft Proof\hspace{2mm} }}
\newcommand{\TF}{\itt}
\newcommand{\bass}{\begin{assertion}\TF\ } 
\newcommand{\eass}{\end{assertion}}
\newcommand{\bcor}{\begin{corollary}\TF\ }
\newcommand{\ecor}{\end{corollary}}
\newcommand{\bdf} {\begin{definition}\rmt\ }
\newcommand{\edf} {\end{definition}} 
\newcommand{\ble} {\begin{lemma}\TF\ }
\newcommand{\ele} {\end{lemma}}
\newcommand{\bpro}{\begin{proposition}\TF\ } 
\newcommand{\epro}{\end{proposition}} 
\newcommand{\brem}{\begin{remark}\rmt\ }
\newcommand{\erem}{\end{remark}} 
\newcommand{\bte} {\begin{theorem}\TF\ }
\newcommand{\ete} {\end{theorem}}
\newcommand{\qed} {\hfill$\msur\Box\msur$} 
\newcommand{\bay}{\begin{array}}
\newcommand{\eay}{\end{array}} 
\newcommand{\bce}{\begin{center}}
\newcommand{\ece}{\end{center}} 
\newcommand{\bde}{\begin{description}}
\newcommand{\ede}{\end{description}}
\newcommand{\ben}{\begin{enumerate}}
\newcommand{\een}{\end{enumerate}}
\newcommand{\bit}{\begin{itemize}}
\newcommand{\eit}{\end{itemize}}
\newcommand{\ZF}  {\bbox{ZF}}
\newcommand{\ZFC} {\bbox{ZFC}}

\newcommand{\DC}  {\bbox{DC}}
\newcommand{\cont}[1]{\rbox{Term}_{#1}}

\newcommand{\We}{{\sf Weak}}
\newcommand{\dom}{\rbox{dom}\,}
\newcommand{\ord}{\rbox{Ord}}
      \newcommand{\Ord}{\ord}

\newcommand{\rL}{\rbox{L}}
\newcommand{\od}{\rbox{OD}}

\newcommand{\rV}{\rbox{V}}
\newcommand{\ima}{{\hbox{\hspace{2pt}\rmt ''}}}
\newcommand{\hc}{\rbox{HC}}
\newcommand{\cuh}[1]{\dd{#1}{\rbox{CUH}}}
\newcommand{\al} {\alpha} 
\newcommand{\ba} {\beta} 
\newcommand{\ga} {\gamma} 
\newcommand{\da} {\delta}
\newcommand{\Da} {\Delta}

\newcommand{\La} {\Lambda} 
\newcommand{\la} {\lambda} 
\newcommand{\sg} {\sigma} 
\newcommand{\Sg} {\Sigma}
 
\newcommand{\vpi}{\varphi} 
\newcommand{\om} {\omega} 

\newcommand{\fs}[2]{{\bsym\Sigma}^{#1}_{#2}}
\newcommand{\fp}[2]{{\bsym\Pi}^{#1}_{#2}}
\newcommand{\fd}[2]{{\bsym\Delta}^{#1}_{#2}}
\newcommand{\iSigma}{{\mathchar"7106}}
\newcommand{\is}[2]{\iSigma^{#1}_{#2}}
\newcommand{\iPi}{{\mathchar"7105}}
\newcommand{\ip}[2]{\iPi^{#1}_{#2}}
\newcommand{\iDelta}{{\mathchar"7101}}
\newcommand{\id}[2]{\iDelta^{#1}_{#2}}
\newcommand{\ish}[1]{\is{\rbox{HC}}{#1}}
\newcommand{\iph}[1]{\ip{\rbox{HC}}{#1}}
\newcommand{\idh}[1]{\id{\rbox{HC}}{#1}}
\newcommand{\fdh}[1]{\fd{\rbox{HC}}{#1}}
\newcommand{\dX}{\mathord{{\rm X}\hspace{-7pt}{\rm X}}}
\newcommand{\dP}{\mathord{{\rm I}\hspace{-2.5pt}{\rm P}}}
\newcommand{\dZ}{\mathord{{\sf Z}\hspace{-4.5pt}{\sf Z}}}
\newcommand{\dF}[1]{{\mathord{{\rm I}\hspace{-2.5pt}{\rm F}}}_{#1}}
\newcommand{\skri}[1]{{\cal #1}}
\newcommand{\cD}{{\skri D}} 
\newcommand{\cN}{{\skri N}}
\newcommand{\cP}{{\skri P}}
\newcommand{\oP}{{\skri P}^{\hbox{\tiny\rm OD}}}
\newcommand{\cX}{{\skri X}}
\newcommand{\sq}  {\subseteq}
\newcommand{\sqq} {\sqsubseteq}
\newcommand{\sneq}{\subsetneqq}
\newcommand{\cj}  {\;\,\&\;\,}

\newcommand{\lra} {\longrightarrow} 
\newcommand{\llra}{\longleftrightarrow} 
\newcommand{\res} {{\mathbin{\hspace{0.15ex}\restriction\hspace{0.15ex}}}}
\newcommand{\we}  {{\mathbin{\hspace{0pt}^\wedge}}}
\newcommand{\<}{\leq}
\def\>{\geq}
\newcommand{\ti}{\times}
\newcommand{\dm}{$$}
\newcommand{\emps}{\emptyset}
\newcommand{\ang} [1]{\langle #1\rangle}
\newcommand{\ans} [1]{\{\hspace{0.2mm}#1\hspace{0.2mm}\}}
\newcommand{\mins}{\hspace{-1pt}-\hspace{-1pt}}
\newcommand{\dd}[1]{$\kern-0.7mm{#1}\kern-1mm$-}
\newcommand{\noi} {\noindent}
\newcommand{\vom} {\vspace{1mm}}
\newcommand{\its} {\vspace{-1mm}}
\newcommand{\itla}[1]{\item\label{#1}}
\def\top{{\cal T}}
\def\C{}
\renewcommand{\C}[1]{{\rbox{C}}_{#1}}
\font\ess=cmss9
\scriptfont\sffam=\ess
\def\E  {\mathbin{\relf{E}}}
\def\Eo {\mathbin{\relf{E}_0}}
\def\nE {\mathbin{\not{\hspace{-2pt}\relf{E}}}}
\def\nEo{\mathbin{\not{\hspace{-2pt}\relf{E}_0}}}
\def\oE {\mathbin{\overline{\relf{E}}}}
\def\noE{\mathbin{\not{\hspace{-2pt}\overline{\relf{E}}}}}
\def\I{}
\def\J{}
\def\R{}
\renewcommand{\I}[1]{\mathbin{\relf{Q}_{#1} }}
\newcommand{\Ip}[1] {\mathbin{\relf{Q}'_{#1}}}
\renewcommand{\J}[1]{\mathbin{\relf{J}_{#1} }}
\renewcommand{\R}[1]{\mathbin{\relf{R}_{#1} }}
\def\X{}
\renewcommand{\X}[1]{{\sf X}_{#1}}
\def\t{}
\renewcommand{\t}{{t}}

\newcommand{\tk}[1]{{\tau[{#1}]}}
\newcommand{\ovu}{{\hat u}}
\newcommand{\ovv}{{\hat v}}
\newcommand{\msur}{\hspace{-1\mathsurround}}

\newcommand{\unf}{{\hat{f}}}
\newcommand{\ung}{{\hat{g}}}
\newcommand{\col}[1]{{#1}^{<\om}}
\begin{document}
\title{On a Glimm -- Effros dichotomy theorem for Souslin  
relations in generic universes}

\author{Vladimir Kanovei
\thanks{Moscow Transport Engineering Institute}
\thanks{{\tt kanovei@nw.math.msu.su} \ and \ 
{\tt kanovei@math.uni-wuppertal.de}
}
\thanks{Partially supported by AMS grant} 
}
\date{01 August 1995} 
\maketitle
\normalsize

\vfill

\begin{abstract}\vspace{2mm} 
\noi
We prove that if every real belongs to a set generic extension of 
the constructible universe then every $\fs11$ equivalence $\E$ on 
reals either admits a $\fd{}1$ reduction to the equality on the 
set $2^{<\om_1}$ of all countable binary sequences, or continuously 
embeds $\Eo,$ the Vitali equivalence.\vom

The proofs are based on a topology generated by $\od$ sets.
\vfill

\bce
{\bft Acknowledgements}
\ece
\noi
The author is in debt to M.~J.~A.~Larijani, the president of 
IPM (Tehran, Iran), for the support during the initial period of 
work on this paper in May and June 1995.
The author is pleased to thank A.~S.~Kechris for the interest 
to this research and G.~Hjorth and A.~W.~Miller for useful 
information on Solovay model and equivalence relations.
\vfill

%\bce
%{\bft Information}
%\ece
%\noi
%At the moment, this is a technical note, not containing an 
%information on the history of the problem, adequate references, etc.
%
\end{abstract}

\vfill

\newpage

\subsection*{Introduction}

This paper presents a proof of the following theorem:

\bte
\label{main}
\footnote
{\rmt\ It follows from an e-mail discussion between G.~Hjorth 
and the author in May -- July 1995 that G.~Hjorth may have proved 
equal or similar theorem independently.}
 \ Let\/ $\E$ be a\/ $\fs11$ 
equivalence on reals. Assume that\its
\bit
\item[$(\dag)$] each real belongs to a ``virtual'' % set 
generic extension~\footnote
{\rmt\ By a generic extension of some $M$ we always mean 
a set generic extension via a forcing notion $P\in M.$ 
Here the extensions could be different for different reals.} 
of the constructible universe~$\rL.$\its
\eit
Then at least one~\footnote
{\rmt\ If all reals are constructible from one of them then 
the statements are compatible.} 
of the following two statements hold$:$\its
\ben
\def\theenumi{{\rmt\hskip2pt(\Roman{enumi})\hskip2pt}}
\def\labelenumi{\theenumi}
\itla{1} \msur 
$\E$ admits a\/ $\fdh1$ reduction~\footnote
{\rmt\ By $\fdh1$ we denote the class of all subsets of $\hc$ 
(the family of all hereditarily countable sets) which are 
$\id{}1$ in $\hc$ by formulas which may contain 
\underline{reals and countable ordinals} as parameters.} 
to the equality on the set\/ $2^{<\om_1}$ of all countable 
binary sequences$.$\its

\itla{2} \msur $\Eo\sqq\E$ continuously$.$
\een
\ete

\subsubsection*{Remarks on the theorem}

By a {\it ``virtual'' generic extension $\rL$\/} we mean a 
set generic extension, say, $\rL[G],$ which is not necessarily 
an inner class in the basic universe $\rV$ (in other words, 
$G\in\rV$ is not assumed).~\footnote
{\rmt\ The assumption that a set $S\sq\Ord$ belongs to a 
``virtual'' set generic extension of $\rL$ can be adequately 
formalized as follows: {\it there exists a Boolean valued extension 
of $\rL[S]$ in which it is true that the universe is a set generic 
extension of the constructible universe\/}, see Lemma~\ref{44} 
below.}

Notice that the assumption $(\dag)$ of the theorem follows e.\ g. 
from the hypothesis that the universe is a set generic extension 
of $\rL.$ In fact the theorem remains true in a weaker assumption 
that each real $x$ belongs to a ``virtual'' generic extension 
of $\rL[z_0]$ for one and the same real $z_0$ which does not depend 
on $x$.

We refer the reader to Harrington, Kechris, and Louveau~\cite{hkl} 
on matters of the early history of ``Glimm -- Effros'' theorems --- 
those of type: {\it each equivalence of certain class either 
admits a reduction to equality or embeds\/ $\Eo$} --- and relevant 
problems in probability and the measure theory. (However 
Section~\ref{ulm} contains the basic notation.) 

The modern history of the topic began in the 
paper \cite{hkl} where it is proved that each Borel equivalence 
on reals either admits a Borel reduction to the equality on reals 
or embeds $\Eo.$ The proof is based on an advanced tool in 
descriptive set theory, the 
{\it Gandy -- Harrington topology\/} on reals, generated by 
$\is11$ sets. 

Hjorth and Kechris~\cite{hk} found that the case of $\fs11$ 
relations is much more complicated. Some examples have shown 
that one cannot find a reasonable ``Glimm -- Effros'' result for 
$\fs11$ relations simply taking a nonBorel reduction in \ref{1} 
or discontinuous embedding in \ref{2}; it seemms that the equality 
on {\it reals\/} rather than countable binary sequences in \ref{1} 
does not match completely the nature of $\fs11$ relations. 

Hjort and Kechris \cite{hk} suggested the adequate approach: one 
has to take $2^{<\om_1}$ as the domain of the equality in 
\ref{1}. (This approach is referred to as the {\it Ulm -- type 
classification\/} in \cite{hk}, in connection with a 
classification theorem of Ulm in algebra.) 
On this way they proved that the dichotomy \ref{1} vs. \ref{2} 
holds for each $\fs11$ equivalence relation on reals, in the 
assumption of the ``sharps'' hypothesis (and the latter can be 
dropped provided the $\fs11$ relation occasionally has only Borel 
equivalence classes). 

Theorem~\ref{main} of this paper establishes the same result (not 
paying attention on the possible compatibility of \ref{1} and 
\ref{2}) in the completely different than sharps assumption: each 
real belongs to a generic extension of $\rL.$ Of course it is 
the principal problem (we may refer to the list of open problems 
in~\cite{hk}) to eliminate the ``forcing'' assumption and prove 
the result in $\ZFC$. 

One faces much more problems in higher projective classes. In 
fact there exists a sort of upper bound for ``Glimm -- Effros'' 
theorems in $\ZFC.$ Indeed, in a nonwellfounded (of ``length'' 
$\om_1\ti\dZ,$ i. e. $\om_1$ successive copies of the integers) 
iterated Sacks extension~\footnote
{\rmt\ See Groszek~\cite{g94} or Kanovei~\cite{k-sacks} on matters 
of nonwellfounded Sacks iterations.} 
of $\rL$ the $\is12$ equivalence
\dm
x\E y\hspace{6mm}\hbox{iff}\hspace{6mm}\rL[x]=\rL[y]
\dm
neither continuously embeds $\Eo$ nor admits a real--ordinal 
definable reduction to the equality on $\cP(\kappa)$ for a 
cardinal $\kappa$. 

Thus the interest can be paid on classes $\fp11,$ $\fd12,$ 
$\fp12.$ One may expect that $\fd12$ relations admit a theorem 
similar to Theorem~\ref{main}.~\footnote
{\rmt\ 
% The author (unpublished) proved that assuming $(\ast)$ and 
% $\om_1^{\rL[z]}$ is countable for each real $x,$ every $\fd12$ 
% relation either admits a $\fdh2$ reduction to the equality on 
% $2^{<\om_1}$ or continuously embeds $\Eo.$ 
G.~Hjorth informed the author that he had partial results in this 
domain.} 

More complicated relations can be investigated in strong 
extensions of $\ZFC$ or in special models. Hjorth~\cite{h-det} 
proved that in the assumption of ${\bf AD}$ and 
$\rV=\rL[{\rm reals}]$ every equivalence on reals either admits 
a reduction (here obviously a real--ordinal definable reduction) 
to the equality on a set $2^\kappa,$ $\kappa\in\Ord,$ or 
continuously embeds $\Eo.$ Kanovei~\cite{k-sm} proved even 
a stronger result (reduction to the equality on $2^{<\om_1}$) in 
Solovay model for $\ZF+\DC$.

\subsubsection*{The organization of the proof}

Theorem~\ref{main} is the main result of this paper. The proof is 
arranged as follows. 

First of all, we shall consider only the case when $\E$ is a 
lightface $\is11$ relation; if in fact $\E$ is $\is11(z)$ in some 
$z\in\cN$ then this $z$ simply enters the reasoning in a uniform 
way, not influenting substantially any of the arguments. 

The splitting point between the statements \ref{1} and \ref{2} 
of Theorem~\ref{main} is determined in Section~\ref{ulm}. 
It occurs that we have \ref{1} in the assumption that 
\bit
\item[$(\ddag)$] each real $x$ belongs to a ``virtual'' \dd\la 
collapsing generic extension of $\rL$ (for some ordinal $\la$) 
in which $\E$ is closed in a 
topolody generated by $\od$ sets on the set ${\cD\cap\We_\la(\rL)}$ 
of all reals \dd\la weak over $\rL.$ (We say that $x\in \cD$ is 
{\it\dd\la weak over\/} $\rL$ iff it belongs to a \dd\al 
collapsing extension of $\rL$ for some $\al<\la$.)
\eit
On the opposite side, we have \ref{2} provided the assumption 
$(\ddag)$ fails. 

Both sides of the proof depend on properties of reals in collapsing 
extensions close to those of Solovay model. The facts we need are 
reviewed in Section~\ref{clos}. 

Section~\ref{ssif} proves assertion \ref{1} of Theorem~\ref{main} 
assuming $(\ddag).$ The principal idea has a semblance of the 
corresponding parts in \cite{hkl} and especially 
\cite{hk}~\footnote
{\rmt\ Yet we use a technique different from the approach of 
\cite{hk}, completely avoiding any use of recursion theory.}
: in the assumption 
of $(\ddag),$ each \dd\la weak over $\rL$ real in the relevant 
``virtual'' \dd\la collapsing extension belongs to a set (one and 
the same 
for all \dd\E equivalent reals) which admits a characterization 
in terms of an element of $2^{<\om_1}.$ An absoluteness argument 
allows to extend this fact to the universe of Theorem~\ref{main}. 

Sections \ref{prod} and \ref{or} prove \ref{2} of 
Theorem~\ref{main} in the assumption that $(\ddag)$ {\it fails\/} 
(but $(\dag)$ still holds, as Theorem~\ref{main} assumes). In fact 
is this case $\E$ is {\it not\/} closed on the set 
$\cD\cap\We_\la(\rL)$ in a ``virtual'' \dd\la collapsing extension 
of $\rL$ for some $\la.$ This suffices to see that $\E$ embeds 
$\Eo$ continuously in the ``virtual'' universe; moreover, $\E$ 
embeds $\Eo$ in a certain special sense which can be expressed by 
a $\is12$ formula (unlike the existence of an embedding in general 
which needs $\is13$). We conclude that $\E$ embeds $\Eo$ in the 
universe of Theorem~\ref{main} as well by Shoenfield.   

The construction of the embedding of $\Eo$ into $\E$ follows the 
principal idea of Harrington, Kechris, and Louveau~\cite{hkl}, yet 
associated with another topology and arranged in a different way. 
(In particular we do not play the strong Choquet game to define 
the necessary sequence of open sets.) \vspace{4mm}

\noi
{\bf Important remark} \\[1mm] 
It will be more convenient to consider $\cD=2^\om,$ the {\em 
Cantor space\/}, rather than $\cN=\om^\om,$ as the basic Polish 
space for which Theorem~\ref{main} is being proved.

\newpage

\subsection{Approach to the proof of the main theorem}
\label{ulm}

First of all, we shall prove only the ``lightface'' case of the 
theorem, so that $\E$ will be supposed to be a $\is11$ equivalence 
on reals. The case when $\E$ is $\is11[z]$ for a real $z$ does not 
differ much: the $z$ uniformly enters the reasoning. 

By ``reals'' we shall understand points of the 
{\it Cantor set\/} $\cD=2^\om$ rather than the {\it Baire space\/} 
$\cN=\om^\om;$ this choice is implied by some technical reasons. 

The purpose of this section is to describe how the two cases of 
Theorem~\ref{main} will appear. This needs to recall some 
definitions. 

\subsubsection{Collapsing extensions}
\label{ce}

Let $\al$ be an ordinal. Then $\col\al$ is the forcing to collapse 
$\al$ down to $\om.$ If $G\sq\col\al$ is \dd{\col\al}generic over 
a transitive model $M$ ($M$ is a set or a class) then $f=\bigcup G$ 
is a function from $\om$ onto $\al,$ so that $\al$ is countable in 
$M[G]=M[f].$ Functions $f:\om\,\lra\,\al$ obtained this way will be 
called \dd{\col\al}{\em generic over\/ $M$.} 

By {\em \dd\la collapse universe hypothesis\/}, \cuh\la{} in 
brief, we shall mean the following assumption: $\rV=\rL[f_0]$ 
for a \dd{\col\la}generic over $\rL$ collapse function 
$f_0\in\la^\om$.

By the assumption of Theorem~\ref{main}, each real $z$ belongs to 
a ``virtual'' \dd{\col\la}generic extension of $\rL,$ the 
constructible universe, for some ordinal $\la.$ Such an extension 
satisfies \cuh\la. 

\brem
\label{rr}
The extension is not necessarily supposed to be an inner class in 
the universe of Theorem~\ref{main}, see Introduction.\qed
\erem
A set is \dd\la{\em weak over $M$} ($\la$ an ordinal in a model 
$M$) iff it belongs to a ``virtual'' \hbox{\dd{\col\al}generic} 
extension of $M$ for some $\al<\la.$ 
We define
\dm
\We_\la(M)=\ans{x:x\,\hbox{ is \dd\la weak over }\,M}\,.
%and ${\weak pX=X\cap\We[p]}$ for every $X$. We put $\We=\We[\emps],$ 
%$\eweak X=\weak\emps X=\ans{x\in X:x\,\hbox{ is weak over }\,\rL}$.
%
\dm
In the assumption \cuh\la{}, reals in $\We_\la(\rL)$ behave 
approximately like all reals in Solovay model. 

\subsubsection{The $\protect\od$ topology}
\label{top}

In $\ZFC,$ Let $\top$ be the topology generated on a given set 
$X$ (for instance, $X=\cD=2^\om,$ the Cantor set) by all $\od$ 
subsets of $X.$ $\top^2$ is the product of two copies of $\top,$ 
a topology on $\cD^2$.

This topology plays the same role in our consideration as the 
Gandy -- Harrington topology in the proof of the classical 
Glimm -- Effros theorem (for Borel relations) in 
Harrington, Kechris, and Louveau~\cite{hkl}. In particular, it 
has similar (although not completely similar: some special 
\dd{\is11}details vanish) properties. 

We define $\oE$ to be the \dd{\top^2}closure of $\E$ in $\cD^2.$ 
Thus ${x\noE y}$ iff there exist $\od$ sets $X$ and $Y$ containing 
resp. $x$ and $y$ and such that ${x'\nE y'}$ for all ${x'\in X,}$ 
${y'\in Y}.$ Obviously $X$ and $Y$ can be chosen as \dd\E 
invariant (simply replace them by their \dd\E saturations), and 
then $Y$ can be replaced by the complement of $X,$ so that 
\dm
x\oE y\;\;\llra\;\;\forall\,X\;
[\,X\hbox{ is }\od \cj X\hbox{ is \dd\E invariant}\;\,
\lra\;\,(x\in X\;\llra\;y\in X)\,]\,.
\dm
Therefore $\oE$ is an $\od$ equivalence on $\cD$. 

\subsubsection{The cases}
\label{cases}

In \cite{hkl}, the two cases are determined by the equality 
$\E=\oE:$ if it holds that $\E$ admits a Borel reduction on 
$\Da(\cD),$ otherwise $\E$ embeds $\Eo.$ Here the splitting 
condition is a little bit more complicated. First of all, we 
have to consider the equality in different universes. Second, 
the essential domain of the equivalence is now a proper subset of 
$\cD,$ the set of all weak reals.\vspace{2mm}

\noi
{\bf Case\ 1.}\ \ For each real $z,$ there exist an ordinal 
$\la$ and a ``virtual'' \dd{\col\la}generic extension $V$ of the 
constructible universe $\rL$ containing $z$ such that the 
following is true in $V:$ $\E$ coincides with $\oE$ on 
$\cD\cap\We_\la(\rL)$ and $x$ is \dd\la weak over $\rL$.\vspace{2mm}

(Notice that, for a $\is11$ binary relation $\E,$ the assertion 
that $\E$ is an equivalence is $\ip12,$ therefore absolute for 
all models with the same ordinals, in particular for $\rL$ and 
all generic extensions of $\rL$.)\vspace{2mm}

\noi
{\bf Case\ 2.}\ \ Not Case 1. 

\bte
\label{mt}
Suppose that each real belongs to a ``virtual'' generic extension 
of\/ $\rL.$ Then, for the given\/ $\is11$ equivalence relation\/ 
$\E,$ we have\/ \its
\bit
\item[--] assertion 
\ref{1} of Theorem~\ref{main} in Case 1,\hfill and\hfill\its
% we have\/ 

\item[--] assertion \ref{2} of Theorem~\ref{main} in Case 2.
\eit
\ete

This is how Theorem~\ref{main} well be proved.

\newpage
\newcommand{\doS}{{\underline S}}
\subsection{On collapsing extensions}
\label{clos}

In this section, we fix a limit constructible cardinal $\la.$ The 
purpose is to establish some properties of \dd\la collapsing 
generic extensions (= the universe under the hypothesis 
\hbox{\cuh\la}). It will be shown that weak ponts (introduced in 
Section~\ref{ulm}) behave approximately like all reals in 
Solovay model. 

\subsubsection{Basic properties}
\label{bp}

We recall that a set $S$ is \dd\la{\em weak over $M$} iff $S$ 
belongs to an \dd{\col\al}generic extension of the model $M$ 
for some $\al<\la$.

The hypothesis \cuh\la{} (the one which postulates that the 
universe is a \dd\la generic extension of $\rL$) will be assumed 
during the reasoning, but we shall not mind to specify \cuh\la{} 
in all formulations of theorems. 

%We recall that 
%${\We_\la(M)=\ans{x:x\,\hbox{ is weak over }\,M}}$ and 
%${\weak pX=X\cap\We_\la[S]}$ for every $X$. We put 
%$\We_\la=\We_\la[\emps],$ 
%$\eweak X=\weak\emps X=\ans{x\in X:x\,\hbox{ is weak over }\,\rL}$.

\bpro
\label{col}
Assume\/ \cuh\la. Let\/ $S\sq \Ord$ be\/ \dd\la weak over\/ $\rL.$ 
Then\its
\ben
\def\theenumi{{\arabic{enumi}}}
\def\labelenumi{{\rmt\theenumi}.}
\itla{sm1}
The universe\/ $\rV$ is a\/ \dd{\col\la}generic extension of 
$\rL[S]$.\its

\itla{sm2} 
If\/ $\Phi$ is a sentence containing only sets in\/ $\rL[S]$ as 
parameters then\/ $\La$ {\rm(}the empty sequence\/{\rm)} decides\/ 
$\Phi$ in the sense of\/ $\col\la$ as a forcing notion over 
$\rL[S]$.\its

\itla{sm4} 
If a set\/ $X\sq\rL[S]$ is\/ $\od[S]$ then\/ $X\in\rL[S]$.
\een 
\epro
($\od[S]=S$\msur{}--{\it ordinal definable\/}, that is, definable 
by an \dd\in formula having $S$ and ordinals as parameters.) 
The proof (a copy of the proof of Theorem 4.1 in Solovay~\cite{sol}) 
is based on several lemmas, including the following crucial lemma:

\ble
\label{44}
Suppose that\/ $P\in\rL$ is a p.o. set, and a set\/ $G\sq P$ is\/ 
\dd Pgeneric over\/ $\rL.$ Let\/ $S\in \rL[G],$ $S\sq\ord.$ 
Then there exists a set $\Sg\sq P,$ $\Sg\in \rL[S]$ such that\/ 
$G\sq\Sg$ and\/ $G$ is\/ \dd{\Sg}generic over\/ $\rL[S]$.
\ele
\proof{}of the lemma. We extract the result from the proof of 
Lemma 4.4 in \cite{sol}. 

{\it We argue in $\rL[S]$}. 

Let $\doS$ be the name for $S$ in the language of the forcing $P$.

Define a sequence of sets 
$A_\al\sq P\;\;(\al\in\Ord)$ by induction on $\al$.\its
\ben
\def\theenumi{{\rmt\hskip1pt(A\arabic{enumi})\hskip1pt}}
\def\labelenumi{\theenumi}
\itla{aa1}\msur
$p\in A_0$ iff either $\sg\in S$ but $p$ forces (in $\rL$ and in 
the sense of $P$ as the notion of forcing) $\sg\not\in\doS,$ or 
$\sg\not\in S$ but $p$ forces $\sg\in\doS$ \ \ --- \ \ for some 
$\sg\in\Ord$.\its

\itla{aa2}\msur
$p\in A_{\al+1}$ iff there exists a dense set $D\sq P,$ $D\in \rL$ 
such that every $q\in D$ satisfying $p\<q$ (means: 
$q$ is stronger than $p$) belongs to $A_\al$.\its

\itla{aa3}
If $\al$ is a limit ordinal then $A_\al=\bigcup_{\ba<\al}A_\ba$.\its
\een

The following properties of these sets are easily verifiable 
(see Solovay \cite{sol}): first, if\/ $p\in A_\al$ and\/ 
$p\< q\in P$ then\/ $q\in A_\al$, second, 
if\/ $\ba<\al$ then\/ $A_\ba\sq A_\al$. 

Since each $A_\al$ is a subset of $P,$ it follows that 
$A_\da= A_{\da+1}$ for some ordinal $\da.$ We put 
$\Sg=P\setminus A_\da.$ Thus $\Sg$ intends to be the set of all 
conditions $p\in P$ which do not force something about $\doS$ 
which contradicts the factual information about $S$. 

We prove, following \cite{sol}, that $\Sg$ is as required. This 
yields a pair of auxiliary facts.\vom

$(\Sg1)$ $G\sq\Sg$.\vom 

\noi 
Indeed assume on the contrary that 
$G\cap A_\ga\not=\emps$ for some $\ga.$ Let $\ga$ be the least 
such an ordinal. Clearly $\ga$ is not limit and $\ga\not=0;$ let 
$\ga=\al+1.$ Let $p\in A_\ga\cap G.$ Since $G$ is generic, 
definition~\ref{aa2} implies $G\cap A_\al\not=\emps,$ 
contradiction.\vom

$(\Sg2)$ {\it If\/ $D\in \rL$ is a dense subset of\/ $P$ then\/ 
$D\cap\Sg$ is a dense subset of\/ $\Sg$}.\vom

\noi
This is easy: if $p\in\Sg$ then $p\not\in A_{\da+1};$ hence 
by~\ref{aa2} there exists $q\in D\setminus A_\da,$ 
%such that 
$q\>p$. 

We prove that $G$ is \dd\Sg generic 
over $\rL[S].$ Let $D\in\rL[S]$ be a dense subset of $\Sg;$ we 
have to check that $D\cap G\not=\emps.$ Suppose that 
$D\cap G=\emps,$ and get a contradiction.

Since ${D\in\rL[S]},$ there exists an \dd\in formula $\Phi(x,y)$ 
containing only ordinals as parameters and such that $\Phi(S,y)$ 
holds in $\rL[S]$ iff $y=D$.

Let $\Psi(G')$ be the conjunction of the following formulas:\its 
\ben
\def\theenumi{\hskip2pt(\arabic{enumi})\hskip2pt}
\def\labelenumi{\theenumi}
\itla{1)}
\msur $S'=\doS[G']$ (the interpretation of the ``term'' $\doS$ via 
$G'$)\hspace{1mm}{} is a set of ordinals, and there exists unique 
$D'\in\rL[S']$ such that $\Phi(S',D')$ holds in $\rL[S']$;\its

\itla{2)}\msur
$D'$ is a dense subset of $\Sg'$ where $\Sg'$ is the set obtained 
by applying our definition of $\Sg$ within $\rL[S']$;\its

\itla{3)}\msur
$D'\cap G'=\emps$.\its
\een 
Then $\Psi(G)$ is true in $\rL[G]$ by our assumptions. Let 
$p\in G$ force $\Psi$ over $\rL.$ Then $p\in\Sg$ by $(\Sg1).$ By 
the density there exists $q\in D$ with $p\< q.$ We can consider 
a \dd\Sg generic over $\rL[S]$ set $G'\sq\Sg$ containing $q.$ 
Then $G'$ is also \dd Pgeneric over $\rL$ by $(\Sg1).$  
We observe that $\doS[G']=S$ because $G'\sq\Sg.$ It follows that 
$D'$ and $\Sg'$ (as is the description of $\Psi$) coinside with 
resp. $D$ and $\Sg.$ In particular $q\in D'\cap G',$ a 
contradiction because $p$ forces \ref{3)}. 
\vspace{3mm}\qed

\proof{}of the proposition. 
{\em Item \ref{sm1}\/}. Lemma~\ref{44} (for $P=\col\la$) 
implies that the universe 
is a \dd\Sg generic extension of $\rL[S]$ for a certain tree 
$\Sg\sq\col\la,$ $\Sg\in\rL[S].$ Notice that $\la$ is a cardinal  
in $\rL[S]$ because $S$ is \dd\al weak over $\rL$ where $\al<\la;$ 
on the other hand, $\la$ is countable in the universe by \cuh\la. 
It follows that there exists a condition $u\in G$ such that the 
set of all \dd\la branching points of $\Sg$ is cofinal over $u$ 
in $\Sg.$ In other words, the set $\ans{v\in\Sg:u\sq v}$ includes 
in $\rL[S]$ a cofinal subset order isomorphic to $\col\la$. 

{\em Items \ref{sm2} and \ref{sm4}\/}. It suffices to 
refer to item \ref{sm1} and argue as in the proofs of Lemma 3.5 and 
Corollary 3.5 in \cite{sol} for $\rL[S]$ as the initial model.\qed

\subsubsection{Coding of reals and sets of reals in the model}
\label{coding}

We let $\dF\al(M)$ be the set of all \dd{\col\al}generic over $M$ 
functions $f\in \al^\om$. 

The following definitions intend to introduce a useful coding 
system for reals (i.\ e. points of $\cD=2^\om$ in this research) 
and sets of reals in collapsing extensions.

Let $\al\in\Ord.$ By $\cont\al$ we denote the set of all indexed 
sets $\t=\ang{\al,\ang{\t_n:n\in\om}}$ -- the ``terms'' -- such that 
${\t_n\sq\col\al}$ for each $n$. 
We put $\cont{<\la}=\bigcup_{\al<\la}\cont\al$. 

``Terms'' $\t\in\cont\al$ are used to code functions 
$C:\al^\om\;\lra\;\cD=2^\om;$ namely, for every $f\in\al^\om$ 
we define $x=\C\t(f)\in\cD$ by: $x(n)=1$ iff $f\res m\in \t_n$ 
for some $m$.

Assume that $\t=\ang{\al,\ang{\t_n:n\in\om}}\in\cont\al,$  
$u\in\col\al,$ $M$ arbitrary. We introduce the sets 
$\X{\t u}(M)=\ans{\C\t(f):u\subset f\in\dF\al(M)}$ and 
$\X\t(M)=\X{\t \La}(M)=\C\t\ima\dF\al(M)$. 
%As above, we let 
%$\X\t[S]=\X\t(\rL[S])$ and $\X\t=\X\t[\emps]=\X\t(\rL);$ the same 
%for $\X{\t u}$.

\bpro
\label{solMb}
Assume\/ \cuh\la. Let $S\sq\Ord$ be\/ \dd\la weak over\/ 
$\rL.$ Then$:$\its
\ben
\def\theenumi{{\arabic{enumi}}}
\def\labelenumi{{\rmt\theenumi}.}
\itla{sm6} 
If\/ $\al<\la,$ $F\sq\dF\al(\rL[S])$ is\/ $\od[S],$ and\/ 
$f\in F,$ then there exists\/ $m\in\om$ such that each\/ 
$f'\in\dF\al(\rL[S])$ satisfying\/ $f'\res m= f\res m$ belongs 
to\/ $F$.\its

\itla{sm5} 
For each real\/ $x\in\cD\cap\We_\la(\rL[S]),$ there exist\/ 
$\al<\la,$ $\t\in\cont\al\cap\rL[S],$ and\/ 
$f\in\dF\al(\rL[S])$ such that\/ $x=\C\t(f)$.
\its

\itla{xl1} 
Each\/ $\od[S]$ set $X\sq\cD\cap\We_\la(\rL[S])$ is a union of 
sets of the form\/ $\X\t(\rL[S]),$ where $\t\in\cont{<\la}\cap\rL[S]$.
\its

\itla{xl2} 
Suppose that\/ ${\t\in\cont\al\cap\rL[S],\;\;\al<\la,}$ and\/ 
${u\in\col\al}.$ Then every\/ $\od[S]$ set\/ 
$X\sq\linebreak[3]{\X{\t u}(\rL[S])}$ is a union of sets of 
the form\/ $\X{\t v}(\rL[S]),$ where\/ $u\sq v\in\col\al$.
\een
\epro
\proof {\em Item \ref{sm6}\/}. We observe that  
$F=\ans{f'\in\al^\om:\Phi(S,f')}$ for an \dd\in formula $\Phi.$ 
Let $\Psi(S,f')$ denote the formula: ``$\La$ \dd{\col\la}forces 
$\Phi(S,f')$ over the universe'', so that 
\dm
F=\ans{f'\in\al^\om:\Psi(S,f')\,\hbox{ is true in }\,\rL[S,f']}.
\dm
by Proposition~\ref{col} (items \ref{sm1} and \ref{sm2}). 
Therefore, since $f\in F\sq\dF\al[S],$ there exists $m\in\om$ such 
that the restriction $u=f\res m$ 
%\in\col\al$ 
\dd{\col\al}forces $\Psi(S,{\hat f})$ over $\rL[S]$ 
where $\hat f$ is the name of the \dd\al collapsing function. 

{\em Item \ref{sm5}\/}. By the choice if $x,$ this real belongs to 
a \dd{\col\al}generic extension of 
$\rL[S].$ Thus $x\in\rL[S,f]$ where $f\in\dF\al(\rL[S]).$ Let 
${\hat x}$ be the name of $x.$ It suffices to define 
$\t_n=\ans{u\in\col\al:u\,\hbox{ forces }\,{\hat x}(n)=1}$ 
and take $\t=\ang{\al,\ang{\t_n:n\in\om}}$.

{\em Item \ref{xl1}\/}. Consider a real $x\in X.$ We use item 2 to 
obtain 
$\al<\la,$ $f\in\dF\al(\rL[S]),$ and\/ $\t\in\cont\al\cap\rL[S]$ 
such that\/ $x=\C\t(f).$ Then we apply item 1 to the $\od[S]$ set 
$F=\ans{f'\in\dF\al[S]:\C\t(f')\in X}$ and the $f$ defined above. 
This results in a condition $u=f\res m\in\col\la$ ($m\in\om$) 
such that $x\in \X{\t u}[S]\sq X.$ Finally the set 
$\X{\t u}[S]$ is equal to $\X{\t'}[S]$ for some other 
$\t'\in \cont\al\cap\rL[S]$.

{\em Item \ref{xl2}\/}. Similar to the previous item.\qed

\newpage

\subsection{The case of closed relations: classifiable points}
\label{ssif}

In this section, we prove the ``case 1'' of Theorem~\ref{mt}. 
Thus let $\E$ be a $\is11$ equivalence relation. 

\subsubsection{Classifiable points}

First of all, we introduce the notion of an \dd\E classifiable 
point.

As usual, $\hc$ denotes the set of all hereditarily countable 
sets. $\ish1$ will denote the collection of all subsets of $\hc$ 
definable in $\hc$ by a parameter-free $\is{}1$ formula. The class 
$\iph1$ is understood the same way, and $\idh1=\ish1\cap\iph1$.

Let us fix a constructible $\idh1$ enumeration 
$\cont{}\cap\rL=\ans{\tk\xi:\xi<\om_1}$ such that each 
$\t\in\cont{}\cap\rL$ has uncountably many numbers 
$\xi<\om_1$ satisfying $\t=\tk\xi.$  
The following lemma gives a more special characterization for 
$\oE,$ the \dd{\top^2}closure of $\E,$ based on this enumeration. 

\ble
\label{har}
Assume\/ \cuh\la. Let\/ $x,\,y\in\cD\cap\We_\la(\rL).$ Then\/ 
${x\oE y}$ if and only if for each\/ $\xi<\om_1$ we have\/ 
$x\in[\X{\tk\xi}(\rL_\xi)]_{\E}\;\llra\;y\in[\X{\tk\xi}(\rL_\xi)]_{\E}$.
\ele
\proof The ``only if'' part follows from the fact that the sets 
$\X{\tk\xi}(\rL_\ga)$ are $\od.$ Let us prove the ``if'' direction. 
Assume that ${x\noE y}.$ There exists an $\od$ set $X$ such that 
$x\in [X]_{\E}$ but $y\not\in [X]_{\E}.$ By 
Proposition~\ref{solMb}, we obtain $x\in \X\t(\rL)\sq [X]_{\E},$ 
where 
$\t=\ang{\al,\,\ang{\t_n:n\in\om}}\in\cont\al\cap\rL,$ $\al<\la.$ 
Since $\la$ is a limit cardinal in $\rL,$ there exists a 
constructible cardinal $\ga,$ $\al<\ga<\la,$ such that 
$\dF\al(\rL)=\dF\al(\rL_\ga).$ Then 
$\t'=\ang{\ga,\,\ang{\t_n:n\in\om}}$ is $\tk\xi$ for some 
$\xi,$ $\ga\<\xi<\om_1.$ 
Then $\X\t(\rL)=\X{\tk\xi}(\rL_\xi)$.\qed
\vspace{4mm} 

For each $x\in\cD,$ we define $\vpi_x\in 2^{\om_1}$ as 
follows: $\vpi_x(\xi)=1$ iff $x\in [\X{\tk\xi}(\rL_\xi)]_{\E}$. 

\bdf
\label{psi}
We introduce the notion of a \hbox{\dd\E c}lassifiable 
point.
We let $T$ be the set of all triples $\ang{x,\psi,t}$ such that 
$x\in\cD,$ $\psi\in 2^{<\om_1},$  
$\t\in\cont\al\cap\rL_{\ga}[\psi],$ where 
$\al<\ga=\dom\psi<\om_1,$ and the following conditions \ref{aa} 
through \ref{cc} are satisfied. %\its 
\ben
\def\theenumi{\hskip2pt{\rmt (\alph{enumi})}\hskip2pt}
\def\labelenumi{\theenumi} 
\itla{aa} \msur 
$\rL_{\ga}[\psi]$ models $\ZFC^-$ (minus the Power Set axiom) so 
that $\psi$ can occur as an extra class parameter in Replacement 
and Separation.\its

\itla{bb}
It is true in $\rL_{\ga}[\psi]$ that $\ang{\La,\La}$ 
forces ${\C\t(\unf)\E\C\t(\ung})$ in the sense of 
$\col\al\!\ti\!\col\al$ as the forcing, where $\unf$ and 
$\ung$ are names for the generic functions in $\al^\om$.\its

\itla{xx}
For each $\xi<\ga,$ $\psi(\xi)=1$ iff 
$x\in [\X{\tk\xi}(\rL_\xi)]_{\E}$ --- so that 
$\psi=\vpi_x\res\ga$.\its

\itla{cc} 
%There exists an \dd{\al^{<\om}}generic over $\rL_{\ga}[\psi]$ 
%function $f\in\al^\om$ such that ${\C\t(f)\E x}$. %\its
\msur
$x$ belongs to $[\X\t(\rL_\ga[\psi])]_{\E}$.
\een
A point $x\in\cD$ is {\em \dd\E classifiable} iff there 
exist $\psi$ and $\t$ such that $\ang{x,\psi,\t}\in T$.\qed
\edf
The author learned from Hjorth and Kechris~\cite{hk} the 
idea of forcing over countable models to get a $\id{}1$ reduction 
function, the key idea of this definition.

\ble
\label{def}
$T_{\E}$ is a\/ $\idh1$ set\/ {\rmt(provided\/ $\E$ is $\is11$)}.
\ele
\proof Notice that conditions \ref{aa} and \ref{bb} in 
Definition~\ref{psi} are $\idh1$ because they reflect truth 
within $\rL_\ga[\psi]$ and the enumeration $\tk\xi$ was chosen 
in $\idh1$. 

Condition \ref{cc} is obviously $\ish1$ (provided $\E$ is at least 
$\is12$), so it remains to convert it also to 
a $\iph1$ form. Notice that in the assumption of \ref{aa} and 
\ref{bb}, the set $X=\X\t(\rL_{\ga}[\psi])$ consists of pairwise 
\dd\E equivalent points. 

(Indeed, consider a pair of \dd{\col\al}generic over 
$\rL_{\ga}[\psi]$ functions $f,\,g\in\al^\om$ (not necessarily a 
{\it generic pair\/}). Let $h\in\al^\om$ be an \dd{\col\al}generic 
over both $\rL_{\ga}[\psi,f]$ and $\rL_{\ga}[\psi,g]$ function. 
Then, by \ref{bb}, ${\C\t(h)\E\C\t(f)}$ holds in 
$\rL_{\ga}[\psi,f,h],$ therefore in the universe by Shoenfield. 
Similarly, ${\C\t(h)\E\C\t(g)}.$ It follows that 
${\C\t(f)\E\C\t(g)},$ as required.) 

Therefore \ref{cc} is equivalent to 
the formula 
$\forall\,y\in \X\t(\rL_{\ga}[\psi])\;(x\E y)$ because 
$\X\t(\rL_{\ga}[\psi])$ is not empty. This is 
clearly $\iph1$ provided $\E$ is at least $\ip12$. 

Let us consider \ref{xx}. The right--hand side of the equivalence 
``iff'' in \ref{xx} is $\is11$ with inserted $\idh1$ functions, 
therefore $\idh1.$ It follows that \ref{xx} itself is 
$\idh1$.~\footnote
{\rmt\ Here we do not see how to weaken the assumption that 
$\E$ is $\is11;$ even if the relation is $\ip11,$ \ref{xx} becomes 
$\idh2$.}
\qed

\subsubsection{The classification theorem}
\label{ct}

The following lemma will allow to define a $\idh1$ reduction of 
the given $\is11$ equivalence relation $\E$ to the equality on 
$2^{<\om_1}$.

\ble
\label{key2}
In the assumption of Case 1 of Subsection~\ref{cases}, each 
point\/ $x\in\cD$ is\/ \dd\E classifiable. 
%moreover, there exist\/ $\psi$ and\/ $\t$ such that\/ 
%$T_{\E}(x,\psi,t)$ and $\dom\psi<\om_1^{\rL[x]}$.
\ele
\proof Let $x\in\cD.$ By the assumption of Case 1, there exist an 
ordinal $\la$ and a ``virtual'' \dd{\col\la}generic extension $V$ 
of the constructible universe $\rL$ containing $x$ such that $\E$ 
coincides with $\oE$ on $\cD\cap\We_\la(\rL)$ in $V$ and $x$ is 
\dd\la weak over $\rL$ in $V$.

Thus we have the two universes, $V$ and the universe of the lemma, 
with one and the same class of ordinals. Since by Lemma~\ref{def} 
``being \dd\E classifiable'' is a $\ish1,$ therefore $\is12$ 
notion, it suffices to prove that $x$ is \dd\E classifiable in the 
``virtual'' universe $V$. 

We observe that \cuh\la{} is true in $V$. 

{\it We argue in $V$}.

Notice that $\vpi=\vpi_x$ is \dd\la weak over $\rL:$ indeed 
$\vpi\in\rL[x]$ by Proposition~\ref{col} since $\vpi$ is $\od[x]$. 
It follows that $[x]_{\E}$ is $\od[\vpi]$ by Lemma~\ref{har}, 
because $\E=\oE$ on $\cD\cap\We_\la(\rL).$ Therefore by 
Proposition~\ref{solMb}, $x\in\X\t(\rL[\vpi])\sq[x]_{\E}$ 
for some $\t\in\cont\al\cap\rL[\vpi],$ $\al<\la$.%\vspace{1mm}

The model $\rL_{\om_1}[\vpi]$ has an elementary submodel 
$\rL_\ga[\psi],$ where $\ga<\om_1$ and $\psi=\vpi\res\ga,$ 
containing $\t$ and $\al.$ We prove that 
$\ang{x,\psi,\t}\in T_{\E}.$ Since conditions \ref{aa} and 
\ref{xx} of Definition~\ref{psi} obviously hold 
for $\rL_\ga[\psi],$ let us check requirements \ref{bb} and 
\ref{cc}.\vspace{1mm}

{\em We check\/ \ref{bb}.} 
Indeed otherwise there exist conditions $u,\,v\in\col\al$ such 
that $\ang{u,v}$ forces ${\C{\t}(\unf)\nE \C{\t}(\ung)}$ in 
$\rL_{\ga}[\psi]$ in the sense of $\col\al\!\ti\!\col\al$ as the 
notion of forcing. Then $\ang{u,v}$ also forces 
${\C{\t}(\unf)\nE \C{\t}(\ung)}$ in $\rL_{\om_1}[\vpi]$. 
Let us consider an 
\dd{\col\al\!\ti\!\col\al}generic over $\rL[\vpi]$ pair 
$\ang{f,g}\in \al^\om\ti\al^\om$ such that $u\subset f$ and 
$v\subset g.$ Then both $y=\C\t(f)$ and $z=\C\t(g)$ belong to 
$\X \t(\rL[\vpi]),$ so ${y\E z}$ because 
$\X\t(\rL[\vpi])\sq [x]_{\E}$.

On the other hand, ${y\E z}$ is {\em false\/} in 
$\rL_{\om_1}[\vpi,f,g],$ that is, in $\rL[\vpi,f,g],$ by the 
forcing property of $\ang{u,v}.$ Therefore we have ${x\nE y}$ 
(in the ``virtual'' universe $V$) by Shoenfield, 
contradiction.\vspace{1mm}

%{\em We check\/ \ref{xx}.} This is clear: indeed, 
%$\psi=\vpi\res\ga$.\vspace{1mm}

{\em We check\/ \ref{cc}.} Take any \dd{\col\al}generic 
over $\rL[\vpi]$ function $f\in\al^\om.$ Then $y=\C\t(f)$ belongs 
to $\X\t(\rL[\vpi]),$ hence ${y\E x}.$ On the other hand, $f$ is 
generic over $\rL_{\ga}[\psi]$.\vspace{1mm}

Thus $\ang{x,\psi,\t}\in T_{\E}.$ This means that 
$x$ is \dd\E classifiable, as required.\qed
%\vspace{4mm}

\bdf
\label{U}
Let $x\in\cD.$ It follows from Lemma~\ref{key2} that there 
exists the least ordinal $\ga=\ga_x<\om_1$ such that 
$T_{\E}(x,\vpi_x\res\ga,\t)$ for some $\t.$ We put 
$\psi_x=\vpi_x\res\ga$ and let $\t_x$ denote the least, in the 
sense of the $\od[\psi_x]$ wellordering of $\rL_{\ga}[\psi_x],$ 
``term'' $\t\in\cont{}[\psi_x]\cap \rL_{\ga}[\psi_x]$ which 
satisfies $T_{\E}(x,\psi_x,\t).$ We put 
$U(x)=\ang{\psi_x,\t_x}$.\qed 
\edf

\ble
\label{inv}
If each\/ $x\in\cD$ is\/ \dd\E classifiable then the map\/ $U$ 
is a\/ $\idh1$ reduction of\/ $\E$ to equality. 
\ele
\proof First of all, $U$ is $\idh1$ by Lemma~\ref{def}. 

If $x\E y$ then $U(x)=U(y)$ because Definition~\ref{psi} 
is \dd\E invariant for $x$. 

Let us prove the converse. Assume that $U(x)=U(y),$ that is, in 
particular, $\psi_x=\psi_y=\psi\in 2^{<\om}$ and 
${\t_x=\t_y=\t\in\cont\al[\psi]\cap\rL_{\ga}[\psi],}$ 
where $\al<\ga=\dom\psi<\om_1$. 

By \ref{cc} we have ${\C\t(f)\E x}$ and ${\C\t(g)\E y}$ for 
some \dd{\col\al}generic over $\rL_{\ga}[\psi]$ functions 
$f,\,g\in\al^\om.$ However $\C\t(f)\E \C\t(g)$ (see the proof 
of Lemma~\ref{def}).\qed

\bcor
\label{case1}
{\rmt[\hspace{1pt}The classification theorem\hspace{1pt}]}\\[1mm]
In the assumption of Case 1 of Subsection~\ref{cases}, $\E$ 
admits a\/ $\idh1$ reduction to the equality on $2^{<\om_1}$. 
\ecor
\proof The range of the function $U$ can be covered by a subset 
$R\sq \hc$ (all pairs $\ang{\psi,\t}$ such that ...) which admits 
a $1-1$ $\idh1$ correspondence with $2^{<\om_1}$.\qed\vspace{4mm}

This completes the proof of the ``case 1'' part of Theorem~\ref{mt}.

\newpage

\subsection{$\protect\od$ forcing}
\label{prod}

This section starts the proof of the ``Case 2'' part of 
Theorem~\ref{mt}. At the beginning, we reduce the problem to a 
more elementary form. 

\subsubsection{Explanation}
\label{expl}

Thus let us suppose that each real $x$ belongs to a ``virtual'' 
generic extension of $\rL,$ but the assumption of Case 1 in 
Subsection~\ref{cases} fails. 

This means the following. There exists a real $z\in\cD$ such that 
for every ordinal $\la$ 
and a ``virtual'' \dd{\col\la}generic extension $V$ of the 
constructible universe $\rL$ containing $z,$ the following is true 
in $V:$ if $z$ is \dd\la weak over $\rL$ then $\E$ 
{\it does not\/} coincide with $\oE$ on $\cD\cap\We_\la(\rL)$.

% Let us fix this real $z$.

We know indeed that 
%by the previous assumption 
$z$ belongs to a 
``virtual'' generic extension of $\rL.$ Therefore there 
exists a limit constructible cardinal $\la$ such that $z$ belongs 
to a \dd{\col\la}generic extension $V$ of $\rL$ and $z$ is weak 
in $V.$ (Simply take $\la$ sufficiently large.) 

Let us fix $\la$ and $V$. 
%
%We have \cuh\la{} in $V.$ Furthermore, it follows from what is 
%said above that $\E\sneq\oE$ on $\cD\cap\We_\la(\rL)$ in $V$.
%
As the condensed matter of this reasoning, we obtain 
\bit
\item\msur
$V$ is a ``virtual'' \dd{\col\la}generic extension of $\rL,$ 
$\la$ is a limit cardinal in $\rL,$ and $\E\sneq\oE$ 
on $\cD\cap\We_\la(\rL)$ in $V$. 
\eit
This is the description of the starting position of the proof of 
the ``Case 2'' part of Theorem~\ref{mt}. The aim is to see that in 
this case $\E$ continuously embeds $\Eo$ in the universe of 
Theorem~\ref{mt}.

The general plan will be first to prove that $\E$ continuously 
embeds $\Eo$ {\it in the auxiliary ``virtual'' universe\/} $V,$ 
and second, to get the result in the universe of Theorem~\ref{mt} 
by Shoenfield. 

After a short examination, one can see a problem in this plan: 
the existence of a continuous embedding $\Eo$ into $\E$ is in fact 
a $\is13$ statement:
\dm
\exists\,\hbox{ continuous }1-1\,\;U:\cD\,\lra\,\cD\;
\forall\,x,\,y\in\cD\;
\left[
\bay{cccl}
x\Eo y & \lra & U(x)\E U(y), & \hbox{ and }\\[2mm]

x\nEo y & \lra & U(x) \nE U(y) & 
\eay
\right]
\dm
The lower implication in the square brackets is $\ip11,$ which 
would match the total $\is12,$ but the upper one is $\is11,$ so 
that the total result is $\is13,$ worse than one needs for 
Shoenfield. 

\subsubsection{Special embeddings and proof of the ``Case 2'' 
part of Theorem~\protect\ref{mt}}
% : an outline}
\label{ne}

To overcome this obstacle, we strengthen the upper implication 
to convert it to a $\ip11$ (actually $\id11$) statement. We recall 
that the $\is11$ set $\E\sq\cD^2$ admits a partition 
$\E=\bigcup_{\al<\om_1}\E_{\al}$ onto Borel sets $\E_\al$ -- the 
{\it constituents\/}, uniquely defined as soon as we have fixed a 
$\ip01$ set $F\sq\cD^2\ti\cN$ which projects onto $\E.$ 

\bdf
\label{nice}
A $1-1$ function $\phi:\cD\,\lra\,\cD$ is a {\it special 
embedding\/} of $\Eo$ into $\E$ iff\its
\ben
\def\theenumi{\rmt({\arabic{enumi}})}
\def\labelenumi{\theenumi}
\itla{cp}
there exists an ordinal $\al<\om_1$ such that 
$\ang{\phi(0^k\we 0\we z),\phi(0^k\we 1\we z)}\in\E_\al$\\ for all 
$z\in\cD$ and $k\in\om$, \hfill and \its

\item for all $x,\,y\in\cD,$ if $x\nEo y$ then $\phi(x)\nE\phi(y)$.
\qed
\een
\edf
%\vspace{1mm}
($0^k$ is the sequence of $k$ zeros.) 
First of all, let us see that a special embedding is an embedding in 
the usual sense. We have to prove that $x\Eo y$ implies 
$\phi(x)\E \phi(y).$ We say that a pair of points $x,\,y\in\cD$ 
is a {\it neighbouring pair\/} iff there exist $k\in\om$ and 
$z\in\cD$ such that $x=0^k\we 0\we z$ and $y=1^k\we 1\we z$ or vice 
versa. Obviously a neighbouring pair is \dd\Eo equivalent. Conversely, 
if $x\Eo y$ then $x$ and $y$ can be connected by a finite chain 
of neighbouring pairs in $\cD.$ Therefore condition~\ref{cp} actually 
suffices to guarantee that $x\Eo y\,\lra\,\phi(x)\E\phi(y)$.

Obviously the existence of a {\it special\/} embedding of $\Eo$ 
into 
$\E$ is a $\is12$ property. Thus, by Shoenfield, to complete the 
proof of the ``Case 2'' part of Theorem~\ref{mt}, it suffices to 
prove the following theorem (and apply it in the auxiliary 
``virtual'' universe $V$).

\bte
\label{mtv}
Assume\/ \cuh\la. 
Let\/ $\E$ be a\/ $\is11$ relation and\/ $\E\sneq\oE$ on\/ 
$\cD\cap\We_\la(\rL).$ Then\/ $\Eo$ admits a special continuous 
embedding into\/ $\E$.
\ete

This theorem is being proved in this and the next section. During 
the course of the proof, we assume \cuh\la{} and fix a $\is11$ 
equivalence $\E$ satisfying $\E\sneq\oE$ on the set 
$\cD\cap\We_\la(\rL)$ 
(although the last assumption will not be used at the beginning). 

In this section, we consider important interactions between $\E$ 
and $\oE.$ The next section defines the required embedding. This 
will complete the proof of theorems~\ref{mtv} and \ref{mt}, and 
Theorem~\ref{main} -- the main theorem.

\subsubsection{$\od$ topology and the forcing}
\label{tforc}

We recall that $\top$ be the topology generated by all $\od$ 
sets.
% As usual $\top=\top[\emptyset]$.

A set $X$ will be called \dd\top{\it separable\/} if the $\od$ 
power set ${\oP(X)=\cP(X)\cap\od}$ has only 
countably many different $\od$ subsets. 

\ble
\label{dizl}
Assume \cuh\la. Let\/ $\al<\la$ and\/ $\t\in\cont\al\cap\rL.$ Each
set\/ $X=\X\t(\rL)$ satisfying\/ $X\sq\cD\cap\We_\la(\rL)$ 
is\/ \dd\top separable.
\ele
\proof By Proposition~\ref{solMb} every $\od$ subset of $X$ is 
uniquely determined by an $\od$ subset of $\col\al.$ Since each 
$\od$ set $S\sq\col\al$ is constructible, we obtain an $\od$ map 
$h:\al^+\,\hbox{ onto }\,\oP(X),$ where $\al^+$ is the least 
cardinal in $\rL$ bigger than $\al.$ Therefore $\oP(X)$ has 
\dd{\<\al^{++}}many $\od$ subsets. It remains to notice that 
$\al^{++}<\la$ because $\la$ is a limit cardinal in $\rL,$ but 
$\la$ is countable in the universe.\qed\vspace{4mm}

Let $\dX=\ans{X\sq\cD:X\,\hbox{ is }\,\od\;\hbox{ and nonempty}\,}$.

Let us consider $\dX$ as a forcing notion (smaller sets are 
stronger conditions) for generic extensions of $\rL.$ Of course 
formally $\dX\not\in\rL,$ but $\dX$ is $\od$ order isomorphic to 
a partially ordered set in $\rL$. (Indeed it is known that there 
exists an $\od$ map $\ell:$ ordinals onto the class of all $\od$ 
sets. Since $\dX$ itself is $\od,$ $\dX$ is a 1--1 image of an 
$\od$ set $\dX'$ of ordinals via $\ell.$ By Proposition~\ref{col} 
both $\dX'$ and the \dd{\ell\hspace{0.5pt}}preimage of the order 
on $\dX$ belong to $\rL$.)
\pagebreak[3]

It also is true that a set $G\sq\dX$ is \dd\dX generic over $\rL$ 
iff it nonempty intersects every dense $\od$ subset of $\dX$. 

\bcor
\label{exis}
Assume \cuh\la. If\/ a set ${X\in\dX}$ satisfies\/ 
$X\sq\cD\cap\We_\la(\rL)$ then there exists a\/ \dd\dX generic 
over\/ $\rL$ set\/ $G\sq\dX$ containing $X$.
\ecor
\proof We can suppose, by Proposition~\ref{solMb}, that 
$X=\X\t(\rL)$ where $\t\in\cont\al\cap\rL$ and $\al<\la.$ Now 
apply Lemma~\ref{dizl}.\qed

\ble 
\label{choq-cor}
Assume \cuh\la. If\/ $G\sq\dX$ is a generic over\/ $\rL$ set 
containing the set\/ $\cD\cap\We_\la(\rL)$ then the intersection\/ 
$\bigcap G$ is a singleton $\ans{a}=\ans{a_G}$.
\ele
\proof Assume that this is not the case. Let $\dX'\in\rL$ be a 
constructible p. o. set order isomorphic $\dX$ via an $\od$ 
function $\ell:\dX'\,\hbox{ onto }\,\dX.$ Then $G'=\ell^{-1}(G)$ is 
\dd{\dX'}generic over $\rL.$ We assert that the statement that 
$\bigcap G$ is not a singleton can be converted to a sentence 
relativized to $\rL[G']$. 

(Indeed, it follows from the reasoning in the proof of 
Lemma~\ref{dizl} that $\rL[G']$ is in fact a \dd Pgeneric 
extension of $\rL$ for a certain set $P\in\rL,$ $P\sq\dX'$ of 
a cardinality $\al<\la$ in $\rL.$ The next \dd\rL cardinal 
$\al^+$ is $<\la$ since $\la$ is a limit cardinal in $\rL.$ 
Therefore $G'$ belongs to a \dd{\col{\al^+}}generic extension of 
$\rL,$ so $G'$ is weak. Then by Proposition \ref{col} the universe 
$\rV=\rL[f_0]$ is a \dd{\col\la}generic extension of $\rL[G'].$ 
This is enough to convert any statement about $G'$ in $\rV$ -- 
like the statement: $\bigcap \ell\ima G'$ is not a singleton -- 
to a sentence relativized to $\rL[G']$.)

Then there exists ${X\in \dX},$ ${X\sq\cD\cap\We_\la(\rL)},$ such 
that $\bigcap G$ is not a singleton for {\em every\/} generic over 
$\rL$ set $G\sq \dX$ containing $X.$ We can assume that 
$X=\X\t(\rL),$ where ${\t\in\cont\al\cap\rL},$ $\al<\la.$ Then $X$ 
is \dd\top separable; let $\ans{\cX_n:n\in\om}$ be an enumeration 
of all $\od$ dense subsets of $\oP(X).$ Using 
Proposition~\ref{col} (item~\ref{sm6}), we obtain an increasing 
\dd{\col\al}generic over $\rL$ sequence 
$u_0\sq u_1\sq u_2\sq...$ of $u_n\in\col\al$ such that 
$X_n=\X{\t\hspace{0.1em}{u_n}}(\rL)\in\cX_n.$ Obviously this 
gives an \dd\dX generic over $\rL$ set $G\sq\dX$ containing $X$ 
and all $X_n$.

Now let $f=\bigcup_{n\in\om}u_n;$ $f\in\al^\om$ and $f$ is 
\dd{\col\al}generic over $\rL.$ Then $x=\C\t(f)\in X_n$ for all 
$n,$ so $x\in\bigcap G.$ Since $\bigcap G$ obviously cannot 
contain more than one point, it is a singleton, so we get a 
contradiction with the choice of $X$.\qed\vspace{4mm}

Reals $a_G$ will be called \dd\od{\it generic over\/} 
$\rL$.%\vspace{3mm}

%It is a separate question, perhaps not an easy one, whether a 
%pair $\ang{a,b}$ of 

\subsubsection{The product forcing}

We recall that $\E$ is assumed to be a $\is11$ equivalence on 
$\cD;$ $\oE$ is the closure of $\E$ in 
the topology $\top^2$ (the product of two copies of $\top$). 

For a set ${P\sq\cD^2,}$ we put 
${\pri P=\ans{x:\exists\,y\;P(x,y)}}$ and 
${\prii P=\ans{y:\exists\,x\;P(x,y)}.}$ Notice that if $P$ is 
$\od,$ so are $\pri P$ and $\prii P$. 

The classical reasoning in Harrington, Kechris, and 
Louveau~\cite{hkl} plays on interactions between $\E$ and $\oE.$ 
In the forcing setting, we have to fix a restriction by $\oE$ 
directly in the definition of the product forcing. Thus we 
consider 
\dm
\dP=\dP(\oE)=\ans{P\sq\oE:P\,\hbox{ is }\od\;\hbox{ and 
nonempty and }\,P=(\pri P\ti \prii P)\cap\oE}
\dm 
as a forcing notion. As above for $\dX,$ the fact that formally 
$\dP$ does not belong to $\rL$ does not cause essential problems.

The following assertion connects $\dP$ and $\dX$.

\vbox{
\bass
\label{proe}
Assume \cuh\la. Then\its
\ben
\def\theenumi{\rmt{\arabic{enumi}}}
\item If\/ $P\in\dP$ then\/ $\pri P$ and\/ $\prii P$ belong to 
$\dX$.\its

\item If\/ $X,\,Y\in\dX$ and\/ $P=(X\ti Y)\cap \oE\not=\emptyset$ 
then $P\in\dP$.\its

\itla{i3} 
If\/ $P\in\dP,$ $X\in\dX,$ $X\sq\pri P,$ then 
there exists\/ $Q\in\dP,$ $Q\sq P,$ such that $X=\pri Q.$  
Similarly for $\prii$.
\een
\eass
}\its\its
\proof Set $Q=\ans{\ang{x,y}\in P:x\in X\cj y\oE x}$ in 
item~\ref{i3}.\qed\vspace{4mm}

A set $P\in\dP$ is {\em \dd\dP separable} if 
the set $\dP_{\sq P}=\ans{Q\in\dP:Q\sq P}$ has only countably many 
different $\od$ subsets. 

\ble 
\label{dizl2}
Assume \cuh\la. Let\/ $\t,\,\t'\in\cont{<\la}\cap\rL.$ Suppose 
that the sets\/ $X=\X\t(\rL)$ and\/ $Y=\X{\t'}(\rL)$ satisfy\/ 
$X\cup Y\sq\cD\cap\We_\la(\rL),$ and finally that\/ 
${P=(X\ti Y)\cap\oE}$ 
is nonempty. Then\/ $P\in\dP$ and\/ $P$ is\/ \dd\dP separable.
\ele
\proof $P\in\dP$ by Assertion~\ref{proe}. A proof of the 
\dd\dP separability can be obtained by a minor modification of 
the proof of Lemma~\ref{dizl}.\qed

\ble
\label{dp2oe}
Assume \cuh\la. Let\/ $G\sq \dP$ be a\/ \dd{\dP}generic over\/ 
$\rL$ set containing\/ ${(\cD\cap\We_\la(\rL))^2\cap\oE}.$ 
Then the intersection\/ $\bigcap G$ contains a single point\/ 
$\ang{a,b}$ where\/ $a$ and\/ $b$ are\/ \dd\od generic 
over $\rL$ and $a\oE b$.
\ele
\proof By Assertion~\ref{proe}, both $G_1=\ans{\pri P:P\in G}$ and 
$G_2=\ans{\pri P:P\in G}$ are \dd\od generic over $\rL$ subsets of 
$\dX,$ so that there exist unique \dd\od generic over $\rL$ points 
$a=a_{G_1}$ and $b=a_{G_2}.$ It remains to show that 
$\ang{a,b}\in\oE$.

Suppose not. There exists an \dd\E invariant $\od$ set 
$A$ such that we have $x\in A$ and $y\in B=\cD\setminus A.$ Then 
$A\in G_1$ and $B\in G_2$ by the genericity. There exists a 
condition $P\in G$ such that $\pri P\sq A$ and $\prii B\sq B,$ 
therefore ${P\sq (A\ti B)\cap\oE=\emptyset},$ which is 
impossible.\qed\vspace{4mm}

Pairs $\ang{a,b}$ as in Lemma~\ref{dp2oe} will be called 
\dd\dP{\it generic\/} and denoted by $\ang{a_G,b_G}$.

For sets $X$ and $Y$ and a binary relation $\R{}\,,$ let us write 
${X\R{}Y}$ if and only if 
$\forall\,x\in X\;\exists\,y\in Y\;(x\R{} y)$ \ and \ 
$\forall\,y\in Y\;\exists\,x\in X\;(x\R{} y)$.

\ble
\label{1for2}
Assume \cuh\la. Let\/ $P_0\in\dP,$ $P_0\sq(\cD\cap\We_\la(\rL))^2,$ 
points\/ $a,\,a'\in X_0=\pri P_0$ be\/ 
\dd\od generic over\/ $\rL,$ and\/ ${a\oE a'.}$ There exists a 
point\/ $b$ such that both\/ $\ang{a,b}$ and $\ang{a',b}$ belong 
to\/ $P_0$ and are\/ \dd{\dP}generic pairs.
% so that in particular\/ $a\oE b$ and\/ $a'\oE b$ by Lemma~\ref{dp2oe}.
\ele
\proof By Lemma \ref{dizl2} and Proposition~\ref{solMb} there 
exists a \dd\dP separable set $P_1\sq P_0$ such that 
$a\in X_1=\pri P_1.$  We put $Y_1=\prii P_1;$ then $X_1\oE Y_1,$ 
and $P_1=(X_1\ti Y_1)\cap\oE$. 

We let $P'=\ans{\ang{x,y}\in P_0:y\in Y_1}.$ Then $P'\in \dP$ 
and $P_1\sq P'\sq P_0.$ Furthermore $a'\in X'=\pri P'.$ (Indeed, 
since ${a\in X_1}$ and ${X_1\oE Y_1},$ there exists $y\in Y_1$ 
such that $a\oE y;$ then $a'\oE y$ as well because $a\oE a',$ 
therefore $\ang{a',y}\in P'$.) By Lemma \ref{dizl2} and 
Proposition~\ref{solMb} 
there exists a \dd\dP separable set $P'_1\sq P'$ such that 
$a'\in X'_1=\pri P'_1.$ Then $Y'_1=\prii P'_1\sq Y_1$. 

It follows from the choice of $P$ and $P'$ that $\dP$ admits only 
countably many different dense $\od$ sets below $P_1$ and below 
$P'_1.$ Let $\ans{\cP_n:n\in\om}$ and $\ans{\cP'_n:n\in\om}$ 
be\pagebreak[3] enumerations of both families of dense 
sets. We define sets $P_n,\,P'_n\in\dP\;\;(n\in\om)$ 
satisfying the following conditions:\its
\ben
\def\theenumi{(\roman{enumi})}
\def\labelenumi{\theenumi}
\itla{i}
$a\in X_n=\pri P_n$ \ and \ $a'\in X'_n=\pri P'_n$;\its

\itla{ii}
$Y'_n=\prii P'_n \sq Y_n=\prii P_n$ \ and \ $Y_{n+1}\sq Y'_n$;\its

\itla{iii}
$P_{n+1}\sq P_n\,,\,$ $P'_{n+1}\sq P'_n\,,\,$ 
$P_n\in \cP_{n-2}\,,\,$ and \ $P'_n\in \cP'_{n-2}$.\its
\een
By \ref{iii} both sequences $\ans{P_n:n\in\om}$ and 
$\ans{P'_n:n\in\om}$ are \dd\dP generic over $\rL,$ so by 
Lemma~\ref{dp2oe} they result in two generic pairs,  
$\ang{a,b}\in P_0$ and $\ang{a',b}\in P_0, $ having the first 
terms equal to $a$ and $a'$ by \ref{i} and second terms equal to 
each other by \ref{ii}. Thus 
it suffices to conduct the construction of $P_n$ and $P'_n$.

The construction goes on by induction on $n$.

Assume that $P_n$ and $P'_n$ have been defined. We define 
$P_{n+1}.$ By~\ref{ii} and Assertion~\ref{proe}, the set 
${P=(X_n\ti Y'_n)\cap\oE\sq P_n}$ belongs to $\dP$ and 
$a\in X=\pri P.$ (Indeed, $\ang{a,b}\in P,$ where $b$ satisfies 
$\ang{a',b}\in P'_n,$ because ${a\oE a'}$.) However $\cP_{n-1}$ 
is dense in $\dP$ below $P\sq P_0;$ therefore 
${\pri \cP_{n-1}=\ans{\pri P':P'\in \cP_{n-1}}}$ is dense in $\dX$ 
below\pagebreak[3] 
$X=\pri P.$ Since $a$ is generic, we have $a\in \pri P'$ for 
some $P'\in \cP_{n-1},$ $P'\sq P.$ It\pagebreak[3] 
remains to put $P_{n+1}=P',$ 
and then $X_{n+1}=\pri P_{n+1}$ and $Y_{n+1}=\prii P_{n+1}$.
\pagebreak[3]

After this, to define $P'_{n+1}$ we let 
$P=(X'_n\ti Y_{n+1})\cap\oE,$ etc.\qed 

%\newpage

\subsubsection{The key set}
\label{second}

We recall that \cuh\la{} is assumed, $\E$ is a $\is11$ 
equivalence on $\cD,$ and $\oE$ is the \hbox{\dd{\top^2}closure} 
of $\E$ in $\cD^2.$ By the assumption of Theorem~\ref{mtv}, 
$\E\sneq\oE$ on $\cD\cap\We_\la(\rL).$ This means that there 
exist \dd\oE classes of elements of $\cD\cap\We_\la(\rL)$ which 
include more than one \dd\E class. We define 
the union of all those \dd\oE classes, 
\dm
H=\ans{x\in\cD\cap\We_\la(\rL):\exists\,y\in
\cD\cap\We_\la(\rL)\;(x\oE y\cj x\nE y)}\,.
%=\ans{x:[x]_{\E}\sneq [x]_{\oE}}
\dm
% (Notice that $\cD\cap\We_\la(\rL)$ is \dd\E invariant by Shoenfield.) 
Obviously $H$ is $\od,$ nonempty, and \dd\E invariant {\em inside\/}
$\cD\cap\We_\la(\rL),$ and moreover $H'=H^2\cap\oE\not=\emptyset,$ 
% in fact both projections of $\cH$ are equal to $H,$ 
so that in particular $H'\in\dP$ by Assertion~\ref{proe}. 

\ble
\label{noE} 
Assume \cuh\la. If\/ $a,b\in H$ and\/ $\ang{a,b}$ is\/ 
\dd\dP generic over $\rL$ then ${a\nE b}\hspace{1.5pt}.$
\ele
\proof Otherwise there exists a set $P\in\dP,$ $P\sq H\ti H$ such 
that $a\E b$ holds for {\it all\/} \dd\dP generic $\ang{a,b}\in P.$ 
We conclude that then $a\oE a'\;\lra\;a\E a'$ for all \dd\od generic 
points $a,\,a'\in X=\pri P;$ indeed, take $b$ such that both 
$\ang{a,b}\in P$ and $\ang{a',b}\in P$ are \dd\dP generic, 
by Lemma~\ref{1for2}. In other words the relations $\E$ and 
$\oE$ coincide on the set 
${Y=\ans{x\in X:x\,\hbox{ is \dd\od generic over }\,\rL}\in\dX.}$ 
($Y$ is nonempty by corollaries \ref{exis} and \ref{choq-cor}.) 

Moreover, $\E$ and $\oE$ coincide on the set 
${Z=[Y]_{\E}\cap\cD\cap\We_\la(\rL).}$ Indeed if $z,\,z'\in Z,$ 
${z\oE z'},$ 
then let ${y,\,y'\in Y}$ satisfy ${z\E y}$ and ${z'\E y'}.$ 
Then ${y\oE y'},$ therefore ${y\E y'},$ which implies $z\E z'.$  

We conclude that $Y\cap H=\emptyset$. 

(Indeed, suppose that $x\in Y\cap H.$ Then by definition there 
exists $y\in\cD\cap\We_\la(\rL)$ 
such that ${x\oE y}$ but ${x\nE y}.$ Then ${y\not\in Z}$ because 
$\E$ and $\oE$ coincide on $Z.$ Thus the pair $\ang{x,y}$ belongs 
to the $\od$ set $P=Y\ti [(\cD\cap\We_\la(\rL))\setminus Z].$ 
Notice that $P$ 
does not intersect $\E$ by definition of $Z.$ Therefore 
$\ang{x,y}$ cannot belong to the closure $\oE$ of $\E,$ 
contradiction.) 

But $\emps\not=Y\sq X\sq H,$ contradiction.\qed\vspace{4mm}

Lemma~\ref{noE} is a counterpart of the proposition in 
Harrington, Kechris, Louveau~\cite{hkl} that $\E\res H$ is 
meager in $\oE\res H.$ But in fact the main content of this 
argument in~\cite{hkl} was implicitly taken by Lemma~\ref{1for2}. 

\ble
\label{E}
Assume \cuh\la. Let\/ $X,\,Y\sq H$ be nonempty\/ $\od$ sets 
and\/ ${X\oE Y}.$ There 
exist nonempty\/ $\od$ sets\/ $X'\sq X$ and\/ $Y'\sq Y$ such 
that\/ $X'\cap Y'=\emptyset$ but still\/ $X'\oE Y'$.
\ele
\proof There exist points $x_0\in X$ and 
$y_0\in Y$ such that $x_0\not= y_0$ but ${x_0\oE y_0}.$ 
(Otherwise $X=Y,$ and $\oE$ is the equality on $X,$ which is 
impossible, see the previous proof.) Let $U$ and $V$ be disjoint 
Baire intervals in $\cD$ containing resp. $x_0$ and $y_0.$ 
The sets $X'= X\cap U \cap [Y\cap V]_{\oE}$ and 
$Y'= Y\cap V \cap [X\cap U]_{\oE}$ are as required.\qed
\newpage

\subsection{Embedding $\protect\Eo$ into $\protect\E$}
\label{or}

In this section we end the proof of Theorem~\ref{mtv}. Thus 
we prove, assuming \cuh\la{} and $\E\sneq\oE$ on 
$\cD\cap\We_\la(\rL),$ that $\E$ embeds $\Eo$ via a 
continuous special (see Definition~\ref{nice}) embedding. 

\subsubsection{The embedding}
\label{embed}

By the assumption the set $H$ of Subsection~\ref{second} is 
nonempty; obviously $H$ is $\od.$ By lemmas \ref{dizl}, 
\ref{dizl2}, and Proposition~\ref{solMb} there exists a nonempty 
\dd\top separable $\od$ set $X_0\sq H$ such that the set 
${P_0=(X_0\ti X_0)\cap\oE}$ belongs to $\dP$ and is \dd\dP 
separable. We observe that 
$\pri P_0=\prii P_0=X_0\sq H\sq\cD\cap\We_\la(\rL)$.  

We define a family of sets $X_u\;\;(u\in 2^{<\om})$ satisfying\its
\ben
\def\theenumi{(\alph{enumi})}
\def\labelenumi{\theenumi}
\itla{a} 
$X_u\sq X_0,$ $X_u$ is nonempty and $\od,$ and $X_{u\we i}\sq X_u,$ 
for all $u$ and $i$.\its
\een
In addition to the sets $X_u,$ we shall define relations 
$\J{uv}\sq\cD^2$ for {\em some} pairs $\ang{u,v},$ to provide 
important interconnections between branches in $2^{<\om}$. 

Let $u,\,v\in 2^n.$ We say that $\ang{u,v}$ is a {\em neighbouring 
pair\/} iff $u=0^k\we 0\we r$ and $v=0^k\we 1\we r$ for some $k<n$ 
($0^k$ is the sequence of $k$ terms each equal to $0$) and some 
$r\in 2^{n-k-1}$ (possibly $k=n-1,$ that is, $r=\La$). 

Thus we define sets $\J{uv}\sq X_u\ti X_v$ for all neighbouring pairs 
$\ang{u,v},$ so that the following requirements \ref{b} and  
\ref{d} will be satisfied.\its
\ben
\def\theenumi{(\alph{enumi})}
\def\labelenumi{\theenumi}
\setcounter{enumi}{1}
\itla{b} \msur 
$\J{uv}$ is $\od,$ $\pri \J{uv}=X_u,$ $\prii \J{uv}=X_v,$ and 
$\J{u\we i\,,\,v\we i}\sq \J{uv}$ for every neighbouring pair 
$\ang{u,v}$ and each $i\in\ans{0,1}$.\its

\itla{d} 
For any $k,$ the set $\J k=\J{0^k\we 0\,,\,0^k\we 1}$ is 
\dd\top separable, and $\J k\sq \E_\al$ for some ordinal 
$\al=\al(k)<\om_1$.\its
\een
Notice that if $\ang{u,v}$ is neighbouring then 
$\ang{u\we i,v\we i}$ is neighbouring, but $\ang{u\we i,v\we j}$ 
is not neighbouring for $i\not=j$ (unless $u=v=0^k$ for some $k$). 
% By $0^k$ we denote the sequence of $k$ zeros.

It follows that $X_u \J{uv} X_v,$ therefore 
$X_u\E X_v,$ for all neighbouring pairs $u,\,v.$~\footnote
{\ We recall that $X\J{}Y$ means that 
$\forall\,x\in X\;\exists\,y\in Y\;(x\J{} y)$ and 
$\forall\,y\in Y\;\exists\,x\in X\;(x\J{} y)$.}

\brem
\label{newrem}
Every pair of $u,\,v\in 2^n$ can be tied in $2^n$ by a 
finite chain of neighbouring pairs. It follows that 
${X_u\E X_v}$ and ${X_u\oE X_v}$ hold for {\em all} pairs 
$u,\,v\in 2^n$.\qed
\erem

Three more requirements will concern genericity. 

Let $\ans{\cX_n:n\in\om}$ be a fixed (not necessarily $\od$) 
enumeration of all dense in $\dX$ below $X_0$ subsets of $\dX.$ 
Let $\ans{\cP_n:n\in\om}$ be a fixed enumeration of all dense in 
$\dP$ below $P_0$ subsets of $\dP.$ It is assumed that 
$\cX_{n+1}\sq\cX_n$ and $\cP_{n+1}\sq\cP_n.$ Note that 
${\cX'=\ans{P\in\dP: P\sq P_0\cj \pri P\cap\prii P=\emptyset}}$  
is dense in $\dP$ below $P_0$ by Lemma~\ref{E}, so we can suppose 
in addition that $\cP_0=\cX'$. 

In general, for any \dd\top separable set $S$ let 
$\ans{\cX_n(S):n\in\om}$ be a fixed enumeration of all dense 
subsets in the algebra $\oP(S)\setminus\ans{\emps}$ such 
that $\cX_{n+1}(S)\sq\cX_n(S)$.

We now formulate the three additional requirements.
\its
\ben
\def\theenumi{({\rmt g}\arabic{enumi})}
\def\labelenumi{\theenumi}
\itla{g1} \msur 
$X_u\in \cX_n$ whenever $u\in 2^n$.\its

\itla{g2}
If $u,\,v\in 2^n$ and $u(n\mins 1)\not=v(n\mins 1)$ (that is, the 
last terms of $u,\,v$ are different), then 
$P_{uv}=(X_u\ti X_v)\cap\oE\in \cP_n$.\its

\itla{g3} 
If $\ang{u,v}=\ang{0^k\we 0\we r,0^k\we 1\we r}\in (2^n)^2$ 
then $\J{uv}\in \cX_n(\J k)$.\its
\een
In particular \ref{g1} implies by Corollary~\ref{choq-cor} that 
for any $a\in 2^\om$ the intersection 
$\bigcap_{n\in\om}X_{a\res n}$ contains a single point, denoted 
by $\phi(a),$ which is \dd\od generic over $\rL,$ and the map 
$\phi$ is continuous in the Polish sense. 

\bass
\label{embe}
Assume\/ \cuh\la. 
$\phi$ is a special continuous 1--1 embedding\/ $\Eo$ to $\E$.
\eass
\proof 
Let us prove that $\phi$ is 1--1. Suppose that 
${a\not=b\in 2^\om.}$ Then ${a(n\mins 1)\not=b(n\mins 1)}$ for 
some $n.$ Let ${u=a\res n},$ ${v=b\res n},$ 
so that we have $x=\phi(a)\in X_u$ and $y=\phi(b)\in X_v.$ But 
then the set ${P=(X_u\ti X_v)\cap \oE}$ belongs to $\cP_n$ by 
\ref{g2}, therefore to $\cP_0.$ This implies 
$X_u\cap X_v=\emptyset$ by definition of $\cP_0,$ 
hence $\phi(a)\not=\phi(b)$ as required.  

Furthermore if $a\nEo b$ (which means that $a(k)\not=b(k)$ for 
infinitely many numbers $k$) then $\ang{\phi(a),\phi(b)}$ is 
\dd\dP generic by \ref{g2}, so $\phi(a)\nE \phi(b)$ by 
Lemma~\ref{noE}.

Let us finally verify that 
$\ang{\phi(0^k\we 0\we c),\phi(0^k\we 1\we c)}\in\E_\al$ for all 
$c\in\cD$ and $k\in \om,$ where $\al=\sup_k\al(k)<\om_1.$ 
The sequence of sets 
$W_m=\J{0^k\we 0\we c\res m\,,\,0^k\we 1\we c\res m}\;\;\,(m\in\om)$ 
is then generic over $\rL$ by \ref{g3} in the sense of the forcing 
$\oP(\J k)\setminus\ans{\emps}$ (we recall that 
$\J k=\J{0^k\we 0\,,\,0^k\we 1}$), which is simply a copy of $\dX,$ 
so that by Corollary~\ref{choq-cor} the intersection of all sets 
$W_m$ is a singleton. Obviously the singleton can be only equal to 
$\ang{\phi(0^k\we 0\we c)\,,\,\phi(0^k\we 1\we c)}.$ We conclude 
that $\phi(0^k\we 0\we c)\E_\al \phi(0^k\we 1\we c),$ as 
required.\qed

\subsubsection{Two preliminary lemmas}

Thus the theorem is reduced to the construction of sets $X_u$ 
and $\J{uv}.$ Before the construction starts, we prove a couple 
of important lemmas.

\ble
\label{impo}
Assume\/ \cuh\la. 
Let\/ $X,\,Y\sq\cD\cap\We_\la(\rL)$ be\/ $\od$ sets such that\/ 
$(X\ti Y)\cap\E$ is nonempty. Then\/ $(X\ti Y)\cap\E$ contains a 
weak over\/ $\rL$ point $\ang{x,y}$.
\ele
\proof 
%It is an interesting question, perhaps not a trivial one, 
%whether a pair of weak points can be not weak. In the case we 
%consider, however it is possible to use the absoluteness of 
%$\is11$ formulas. 
First of all, by Proposition~\ref{solMb} we can 
assume that $X=\X{\t}(\rL)$ and $Y=\X{\t'}(\rL),$ where $\t$ and 
$t'$ belong to some $\cont\al\cap\rL,$ $\al<\la.$ Then, since 
$\la$ is a limit \dd\rL cardinal, we have $X=\X{\t}(\rL_\ba)$ and 
$Y=\X{\t'}(\rL_\ba)$ for a suitable $\ba,$ $\al\<\ba<\la.$ Take 
an arbitrary \dd{\col\ba}generic over $\rL$ function 
$f\in\ba^\om.$ Then the statement ${(X\ti Y)\cap\E\not=\emps}$ turns 
out to be a $\is11$ formula with reals in $\rL[f]$ (those coding 
$f,\;\t,\;\t'$) as parameters. 
Notice that all sets in $\rL[f]$ are weak over $\rL,$ so it 
remains to apply Shoenfield.\qed

\ble
\label{comb}
Assume\/ \cuh\la. 
Let\/ $n\in\om,$ and\/ $X_u$ be a nonempty\/ $\od$ set for each\/ 
$u\in 2^n.$ Assume that an\/ $\od$ set\/ $\J{uv}\sq \cN^2$ is 
given for every neighbouring pair of\/ $u,\,v\in 2^n$ so that 
$X_u \J{uv} X_v$.\its
\ben
\def\theenumi{{\rmt\arabic{enumi}.}}
\def\labelenumi{\theenumi}
\item If\/ $u_0\in 2^n$ and\/ $X'\sq X_{u_0}$ is\/ $\od$ and 
nonempty then there exists a system of\/ $\od$ nonempty sets\/ 
$Y_u\sq X_u\;\;(u\in 2^n)$ such that\/ $Y_u \J{uv} Y_v$ holds for 
all neighbouring pairs\/ $u,\,v,$ and in addition $Y_{u_0}=X'$.\its

\item
%Let\/ $u_0,\,v_0\in 2^n$ be a neighbouring pair. 
Suppose that\/ $u_0,\,v_0\in 2^n$ is a neighbouring pair and 
nonempty\/ $\hspace{-1pt}\od\hspace{-1pt}$ 
sets\/ ${X'\sq X_{u_0}}$ and $X''\sq X_{v_0}$ satisfy\/ 
$X' \J{u_0v_0} X''.$ Then there exists a system of\/ $\od$ 
nonempty sets\/ ${Y_u\sq X_u}$ $(u\in 2^n)$ such that\/ 
${Y_u \J{uv} Y_v}$ holds for all neighbouring pairs\/ $u,v,$ and 
in addition\/ $Y_{u_0}=X',\,\;Y_{v_0}=X''.$
\een
\ele
\proof Notice that 1 follows from 2. Indeed take arbitrary $v_0$ 
such that either $\ang{u_0,v_0}$ or $\ang{v_0,u_0}$ is neighbouring, 
and put respectively 
${X''=\ans{y\in X_{v_0}: \exists\,x\in X'\;(x \J{u_0v_0} y)}},$ or 
${X''=\ans{y\in X_{v_0}: \exists\,x\in X'\;(y \J{v_0u_0} x)}}$. 

To prove item 2, we use induction on $n.$ 

For $n=1$ --- then $u_0=\ang{0}$ and $v_0=\ang{1}$ --- 
we take $Y_{u_0}=Y'$ and $Y_{v_0}=Y''$.

The step. We prove the lemma for $n+1$ provided it has been proved 
for $n;\,\,n\>1.$ The principal idea is to divide $2^{n+1}$ on two 
copies of $2^n,$ minimally connected by neighbouring pairs, and handle 
them more or less separately using the induction hypothesis. The 
two ``copies'' are $U_0=\ans{s\we 0:s\in 2^n}$ and 
$U_1=\ans{s\we 1:s\in 2^n}$. 

The only neighbouring pair that connects $U_0$ and $U_1$ is the pair 
of $\ovu=0^n\we 0$ and $\ovv=0^n\we 1.$ If in fact $u_0=\ovu$ 
and $v_0=\ovv$ then we apply the induction hypothesis (item~1) 
independently for the families ${\ans{X_u:u\in U_0}}$ and 
${\ans{X_u:u\in U_1}}$ and the given sets ${X'\sq X_{u_0}}$ and 
${X''\sq X_{v_0}.}$ Assembling the results, we get nonempty $\od$ 
sets ${Y_u\sq X_u\,\;(u\in 2^{n+1})}$ such that ${Y_u \J{uv} Y_v}$ 
for all neighbouring pairs\/ $u,\,v,$ perhaps with the exception of the 
pair of $u=u_0=\ovu$ and $v=v_0=\ovv,$ and in addition 
$Y_{u_0}=X'$ and $Y_{v_0}=X''.$ 
Thus finally $Y_\ovu \J{\ovu\ovv}Y_\ovv$ by the choice of $X'$ and 
$Y'$. 

It remains to consider the case when both $u_0$ and $v_0$ belong 
to one and the same domain, say to $U_0.$ Then we first apply the 
induction hypothesis (item 2) to the family 
${\ans{X_u:u\in U_0}}$ and the sets ${X'\sq X_{u_0}}$ and 
${X''\sq X_{v_0}.}$ This results in a system of nonempty $\od$ 
sets ${Y_u\sq X_u\;\,(u\in U_0);}$ in particular we get an $\od$ 
nonempty set $Y_\ovu\sq X_\ovu.$ We put 
${Y_\ovv=\ans{y\in X_\ovv:\exists\,x\in Y_\ovu\,
(x\J{\ovu\ovv}y)},}$ so that ${Y_\ovu\J{\ovu\ovv}Y_\ovv,}$ 
and apply the induction hypothesis (item 1) to the family 
${\ans{X_u:u\in U_1}}$ and the set $Y_\ovv\sq X_\ovv$.\qed

\subsubsection{The construction}

We put $X_\La=X_0.$ 

Now assume that the sets $X_s\,\;(s\in 2^n)$ 
and relations $\J{st}$ for all neighbouring pairs of $s,\,t\in 2^{\<n}$ 
have been defined, and expand the construction at level $n+1.$ 

We first put $A_{s\we i}=X_s$ for all $s\in 2^n$ and 
$i\in\ans{0,1}.$ We also define $\I{uv}=\J{st}$ for any neighbouring 
pair of $u=s\we i,\,\,v=t\we i$ in $2^{n+1}$ other than the pair 
$\ovu=0^n\we 0,\,\,\ovv=0^n\we 1.$ For the latter one (notice 
that $A_{\ovu}=A_{\ovv}=X_{0^n}$) we put $\I{\ovu\ovv}=\oE,$ 
so that $A_u\I{uv} A_v$ holds for all neighbouring pairs of 
$u,\,v\in 2^{n+1}$ including the pair $\ang{\ovu,\ovv}$.

The sets $A_u$ and relations $\I{uv}$ will be reduced in several 
steps to meet requirements \ref{a}, \ref{b}, \ref{d} and 
\ref{g1}, \ref{g2}, \ref{g3} of Subsection~\ref{embed}.\vom

{\em Part 1}. After $2^{n+1}$ steps of the procedure of 
Lemma~\ref{comb} (item 1) we obtain a system of nonempty $\od$ 
sets $B_u\sq A_u\;\,(u\in 2^{n+1})$ such that still 
$B_u\I{uv} B_v$ for all neighbouring pairs $u,\,v$ in $2^{n+1},$ but 
$B_u\in \cX_{n+1}$ for all $u.$ Thus \ref{g1} is fixed.\vom

{\em Part 2}. To fix \ref{g2}, consider an arbitrary pair of 
$u_0=s_0\we 0,$ $v_0=t_0\we 1,$ where $s_0,\,t_0\in 2^n.$ By 
Remark~\ref{newrem} and density of the set $\cP_{n+1}$ there 
exist nonempty $\od$ sets $B'\sq B_{u_0}$ and $B''\sq B_{v_0}$ 
such that ${P=(B'\ti B'')\cap\oE\in \cP_{n+1}}$ and 
$\pri P=B',$ $\prii P=B'',$ so in particular ${B'\oE B''}.$ Now we 
apply Lemma~\ref{comb} (item 1) separately for the two systems 
of sets, 
${\ans{B_{s\we 0}:s\in 2^n}}$ and ${\ans{B_{t\we 1}:t\in 2^n}}$ 
(compare with the proof of Lemma~\ref{comb}~!), and the sets 
$B'\sq B_{s_0\we 0},$ $B''\sq B_{t_0\we 1}$ respectively. 
This results in a system of nonempty 
$\od$ sets ${B'_u\sq B_u}$ ${(u\in 2^{n+1})}$ satisfying 
${B'_{u_0}=B'}$ and ${B'_{v_0}=B'',}$ so that we have 
${(B'_{u_0}\ti B'_{v_0})\cap\oE\in \cP_{n+1},}$ and still 
$B'_u\I{uv} B'_v$ for all neighbouring pairs $u,\,v\in 2^{n+1},$ 
perhaps with the exception of the pair of 
$\ovu=0^n\we 0,\,\,\ovv=0^n\we 1,$ which is the only one that 
connects the two domains. To handle this exceptional pair, 
note that ${B'_{\ovu} \oE B'_{u_0}}$ and ${B'_{\ovv} \oE B'_{v_0}}$ 
(Remark~\ref{newrem} is applied to each of the two domains), 
so that ${B'_\ovu\oE B'_\ovv}$ since ${B'\oE B''}.$ Finally 
we observe that $\I{\ovu\ovv}$ is so far equal to $\oE$. 

After $2^{n+1}$ steps (the number of pairs $u_0,\,v_0$ to be 
considered) we get a system of nonempty $\od$ sets 
$C_u\sq B_u\;\,(u\in 2^{n+1})$ such that  
$(C_u\ti C_v)\cap\oE\in \cP_{n+1}$ whenever $u(n)\not=v(n),$ 
and still $C_u\I{uv} C_v$ for all neighbouring pairs 
$u,\,v\in 2^{n+1}.$ Thus \ref{g2} is fixed.\vom

{\em Part 3}. We fix \ref{d} for the exceptional neighbouring 
pair of ${\ovu=0^n\we 0},$ ${\ovv=0^n\we 1}.$ Since $\E$ is 
\dd{\top^2}dense in $\oE,$ and ${C_\ovu\oE C_\ovv,}$ the set
${\R{}=(C_\ovu \ti C_\ovv)\cap\E}$ is nonempty. We observe that 
the $\od$ set 
\dm
{\R{}}'=\ans{\ang{x,y}\in {\R{}}:\ang{x,y}\,\hbox{ is weak over }
\,\rL}
\dm
is nonempty, too, by Lemma~\ref{impo}. Then, since 
${\R{}}'\sq\R{}\sq\E,$ the intersection 
${\R{}}''=\R{}'\cap\E_\al$ is nonempty for some $\al<\om_1.$ 
($\E_\al$ is the \dd\al th constituent of the \dd{\is11}set 
$\E$.)  Finally some nonempty $\od$ 
set $\I{}\sq {\R{}}''$ is \dd\top separable by Lemma~\ref{dizl}. 
Consider the $\od$ sets $C'=\pri \I{}\,\,(\sq C_\ovu)$ and 
$C''=\prii \I{}\,\,(\sq C_\ovv);$ obviously $C'\I{} C'',$ so that
$C'\I{\ovu\ovv} C''.$ (We recall that at the moment 
$\I{\ovu\ovv}=\oE.$) Using Lemma~\ref{comb} (item 2) again, we 
obtain a system of nonempty $\od$ sets 
$Y_u\sq C_u\;\,(u\in 2^{n+1})$ such that still $Y_u\I{uv} Y_v$ 
for all neighbouring pairs $u,\,v$ in $2^{n+1},$ and $Y_\ovu=C',$  
$Y_\ovv=C''.$ We re--define $\I{\ovu\ovv}$ by $\I{\ovu\ovv}=\I{}$ 
(then $\I{\ovu\ovv}\sq \E_\al$), 
but this keeps $Y_\ovu\I{\ovu\ovv} Y_\ovv$.\vom

{\em Part 4}. We fix \ref{g3}. Consider a neighbouring pair 
$u_0,\,v_0$ in $2^{n+1}.$ 
Then we have $u_0=0^k\we 0\we r,$ $v_0=0^k\we 1\we r$ for 
some ${k\<n}$ and ${r\in 2^{n-k}}.$ It follows that 
${\Ip{}=\I{u_0v_0}\cap(Y_{u_0}\ti Y_{v_0})}$ is a nonempty 
(since ${Y_{u_0}\I{u_0v_0} Y_{v_0}}$) $\od$ subset of 
$\J k=\J{0^k\we 0\,,\,0^k\we 1}$ by the construction. Let 
$\I{}\sq \Ip{}$ be a nonempty $\od$ set in $\cX_{n+1}(\J k).$ We 
now define $Y'=\pri\I{}$ and $Y''=\prii\I{}$ (then ${Y'\I{}Y''}$ 
and ${Y'\I{u_0v_0}Y''}$) and run Lemma~\ref{comb} (item 2) for the 
system of sets $Y_u\;\,(u\in2^{n+1})$ and the sets 
${Y'\sq Y_{u_0}},$ ${Y''\sq Y_{v_0}}$. After this define 
the ``new'' $\I{u_0v_0}$ by $\I{u_0v_0}=\I{}$.  

Do this consequtively for all neighbouring pairs; the finally 
obtained sets -- let them be $X_u\,\;(u\in 2^{n+1})$ -- are as 
required. The final relations $\J{uv}\;\,(u,\,v\in 2^{n+1})$ can 
be obtained as the restrictions of sets $\I{uv}$ to 
$X_u\ti X_v$.\vom

This ends the construction.
\vspace{4mm}

This also ends the proof of theorems \ref{mtv} and \ref{mt}, 
and Theorem~\ref{main} (the main theorem), see 
Subsection~\ref{ne}.\qed

\end{document}